\documentclass{article}

\usepackage{tikz}
\usepackage{amsmath}
\usepackage{amsfonts}
\usepackage{amssymb}
\usepackage{xspace}
\usepackage{graphicx}
\usepackage[title]{appendix}
\usepackage{color}

\numberwithin{equation}{section}
\makeatletter
\newenvironment{prooof}{\par\noindent{\sc Proof:}
}{\hfill\llap{$\Box$}\vspace{1\baselineskip}\par\noindent}

\newenvironment{proofof}{\par\noindent{\sc Proof}
}{\hfill\llap{$\Box$}\vspace{1\baselineskip}\par\noindent}

\newcommand{\adl}{\vspace{1\baselineskip}}

\newtheorem{theorem}{Theorem}[section]
\newtheorem{proposition}[theorem]{Proposition}
\newtheorem{lemma}[theorem]{Lemma}
\newtheorem{corollary}[theorem]{Corollary}
\newtheorem{remark}[theorem]{Remark}
\newtheorem{definition}[theorem]{Definition}
\newtheorem{example}[theorem]{Example}

\newcommand{\beq}{\begin{equation}}
\newcommand{\eeq}{\end{equation}}

\newcommand{\ba}{\begin{array}}
\newcommand{\ea}{\end{array}}
\newcommand{\bt}{\begin{theorem}}
\newcommand{\et}{\end{theorem}}
\newcommand{\bp}{\begin{proposition}}
\newcommand{\ep}{\end{proposition}}
\newcommand{\bl}{\begin{lemma}}
\newcommand{\el}{\end{lemma}}
\newcommand{\bc}{\begin{corollary}}
\newcommand{\ec}{\end{corollary}}
\newcommand{\bi}{\begin{itemize}}
\newcommand{\ei}{\end{itemize}}
\newcommand{\ben}{\begin{enumerate}}
\newcommand{\een}{\end{enumerate}}
\newcommand{\bpf}{\begin{prooof}}
\newcommand{\epf}{\end{prooof}}
\newcommand{\bpff}{\begin{proofof}}
\newcommand{\epff}{\end{proofof}}
\newcommand{\bdf}{\begin{definition}\rm}
\newcommand{\edf}{\end{definition}}
\newcommand{\br}{\begin{remark}\rm}
\newcommand{\er}{\end{remark}}
\newcommand{\bex}{\begin{example}\rm}
\newcommand{\eex}{\end{example}}
\def\pri{\hbox to 10pt{\hfil\hbox to 0.4pt{\vrule height5pt width0.4pt
                 depth0pt}\vrule width5pt height0.4pt depth0pt\hfil}}
\newcommand{\TC}{{\rm TC}}

\newcommand{\TAT}{{\rm TAT}}

\newcommand{\Lin}{{\rm Lin}}

\newcommand{\gk}{{\bf{k}}}

\newcommand{\gP}{{\bf p}}

\newcommand{\calL}{{\cal L}}

\newcommand{\gR}{{\mathbb R}}
\newcommand{\bI}{{\mathbb I}}
\newcommand{\bJ}{{\mathbb J}}
\newcommand{\bB}{{\mathbb B}}
\newcommand{\bG}{{\mathbb G}}

\newcommand{\Nat}{{\mathbb N}}

\newcommand{\RN}{{\mathbb R}^{N}}
\newcommand{\Rm}{{\mathbb R}^{m}}

\newcommand{\M}{{\cal M}}

\newcommand{\BV}{\mathop{\rm BV}\nolimits}

\def\dist{\mathop{\rm dist}\nolimits}

\def\det{\mathop{\rm det}\nolimits}

\newcommand{\wc}{\rightharpoonup}
\newcommand{\tr}{\mbox{\rm tr}}

\newcommand{\Sph}{{{\mathbb S}}}

\newcommand{\gp}{\mathbf{p}}

\def\mesh{\mathop{\rm mesh}\nolimits}
\newcommand{\ttt}{{\mathfrak t}}

\newcommand{\bbb}{{\mathfrak b}}

\newcommand{\gt}{{\bf t}}
\newcommand{\gn}{{\bf n}}
\newcommand{\gb}{{\bf b}}

\newcommand{\gv}{{\bf v}}

\newcommand{\gr}{{\bf r}}
\newcommand{\gc}{{\bf c}}

\newcommand{\Var}{\mathop{\rm Var}\nolimits}

\def\dd{\mathrm{d}}

\let\a=\alpha
\let\be=\beta

\let\d=\delta
\let\e=\varepsilon
\let\f=\phi
\let\vf=\varphi
\let\g=\gamma

\let\l=\lambda
\let\m=\mu

\let\p=\pi

\let\r=\rho
\let\s=\sigma
\let\t=\theta

\let\vf=\varphi
\let\DD=\Delta

\let\O=\Omega

\let\wih=\widehat

\let\wid=\widetilde

\let\sb=\subset

\let\Lra=\Longrightarrow

\let\fa=\forall
\let\tim=\times

\let\sm=\setminus
\let\ol=\overline

\let\ds=\displaystyle

\let\i=\infty
\let\lm=\limits

\title{\Large \bf The fundamental theorem of the local theory \\ for non-smooth curves}
\author{\it Domenico Mucci and Alberto Saracco
\footnote{%
{\sc Dipartimento di Scienze Matematiche,
Fisiche ed Informatiche (SMFI), Universit\`{a} di Parma,
Parco Area delle Scienze 53/A, I-43124 Parma, Italy.
E-mail: domenico.mucci@unipr.it, alberto.saracco@unipr.it}
}
}
\begin{document}
\date{}
\maketitle
\topskip=1.5truecm
%
{\small{\bf Abstract.} We extend the classical fundamental theorem of the local theory of smooth curves to a wider class of non-smooth data. Curvature and torsion are prescribed in terms of the distributional derivative measures of two given functions of bounded variation. The essentially unique non-smooth curve solution has both finite total curvature and total absolute torsion. In case of continuous data, we preliminarly discuss a more general problem involving a linear system of distributional derivative equations.}
\adl\par\noindent
\adl\par\noindent
\bigskip
\par\noindent
The classical fundamental theorem of the local theory of curves states that given two smooth functions $k(s)$ and $\tau(s)$ on a closed and bounded interval, with $k(s)$ always positive, there exists an essentially unique curve $\gc$ in the Euclidean space $\gR^3$ with curvature and torsion equal to $k$ and $\tau$. 

More precisely, if $I_L=(0,L)$, where $L>0$, and $k(s),\,\tau(s)\in C^1(\ol I_L)$, denoting by $({\bf e}_1,{\bf e}_2,{\bf e}_3)$ the (column) vectors of the canonical basis in $\gR^3$,
it turns out that the linear system of O.D.E.'s
\beq\label{FS-intro}
\left\{ \begin{array}{l}
\gt'(s) =k(s)\,\gn(s) \\
\gn'(s)=-k(s)\,\gt(s)+\tau(s)\,\gb(s) \\
\gb'(s)=-\tau(s)\,\gn(s) \\
\gt(0)={\bf e}_1 \,,\quad \gn(0)={\bf e}_2 \,,\quad \gb(0)={\bf e}_3
\end{array} \right.
\eeq
has a unique solution, and the corresponding matrix valued function 
\beq\label{G-intro}
G(s)=\left(\gt(s) \mid\gn(s)\mid\gb(s)\right) 
\eeq
takes values in the special orthogonal group $SO(3)$ for every $s\in \ol I_L$.
Furthermore, if in addition $k(s)>0$ for every $s\in \ol I_L$, the parametric curve
\[
\gc(s)=0_{\gR^3}+\int_0^s \gt(\s)\,\dd\s \,,\quad s\in\ol I_L
\]
has length $L$ and Frenet frame identified by $G(s)$ for every $s$. Therefore, the solution curve has both finite {\em total curvature} and {\em total absolute torsion} given by the integrals
\[ 
\TC(\gc)=\int_{I_L}k(s)\,ds<\i\,,\quad \TAT(\gc)=\int_{I_L}|\tau(s)|\,ds<\i\,.
\]

The curve is uniquely determined up to rigid motions. However, if we remove the positivity condition and assume that $k(s)=0$ on some non trivial sub-interval of $I_L$, the curve $\gc$ fails to be essentially unique.

In this paper, we wish to extend the previous result to a more general class of given data.
More precisely, referring to Sec. \ref{Sec:notation} for the notation adopted here, we assume that 
\[
\t: I_L\to \gR\,, \quad \f: I_L\to\gR
\]
are two given {\em functions of bounded variation}, and we wish to replace the smooth data $k$ and $\tau$ with the finite Borel measures given by the 
{\em distributional derivatives} $D\t$ and $D\f$. If $\t$ and $\f$ are smooth, we have
\[ D\t=k\,\calL^1\,,\quad D\f=\tau\,\calL^1\,,
\]
where $\calL^1$ is the Lebesgue measure. For that reason, we require in addition that $\t$ is {\em strictly increasing} in $I_L$. 

With these assumptions, the Frenet-Serret system has to be replaced with a linear system involving distributional derivatives, and possible solutions are expected to have a weaker regularity. Therefore, it is reasonable to deal with matrix-valued functions $s\mapsto G(s)\in SO(3)$ whose entries are functions of bounded variation. 

However, simple examples in the planar case (i.e., assuming $\f\equiv 0$) show that if the increasing angle function $\t$ has at least one discontinuity point, then the expected solution $G(s)$ fails to be continuous, and the solution curve $\gc$ has a corner point in correspondence of every jump of $\t$, that produces a jump of the {\em tangent indicatrix} function $\gt(s)$. In addition, it is not immediate how to define the system of distributional derivative equations at the Jump points of the data $\t$ and $\f$.

For these reasons, we shall first analyze the case when both the BV functions $\t$ and $\f$ have a {\em continuous representative}, where we expect to find continuous solutions $s\mapsto G(s)\in SO(3)$ of bounded variation.

After collecting in Sec. \ref{Sec:notation} some background material, in Sec. \ref {Sec:cont} we consider a more general problem with continuous data. More precisely, in any dimension $N\geq 2$, we consider a continuous function of bounded variation $\O : I_L\to Sk(N)$ taking values in the set of skew-symmetric and real valued $N\times N$ matrices $Sk(N)$.
For a given initial condition $\bG$ in the special orthogonal group $SO(N)$, we consider the system
\beq\label{eqND-intro}
\left\{ \begin{array}{l}
DG= - (G^\top)^{-1}D\O \\
G(0)=\bG
\end{array} \right.
\eeq
and we look for solutions $G:I_L\to \gR^{N\tim N}$ such that each entry $G(s)^i_j$ is a function of bounded variation. We are thus implicitly requiring that $G(s)$ is an invertible matrix for every $s$. 

Our first main result is the following
\bt\label{Tintr} The system \eqref{eqND-intro} admits a unique solution, that is given by a continuous function of bounded variation $G:I_L\to SO(N)$.
\et

We did not find in the literature a proof of the latter result, and we refer to \cite{Stallard} for a related problem concerning solutions to differential systems in classes of functions of bounded variation.
More recently, in the Sobolev setting geometric invariants of non-smooth framed curves in the sense of Bishop \cite{Bi} have been analyzed in \cite{BLM}. The papers \cite{FHMP,FHMP2,GF,GLF,GMSV,S} contain applications of the theory of framed curve to ribbons, Cosserat rods, elastic rods. Concerning rod deformations, we refer to \cite{BCHK} for the analysis of the non-interpenetration property, and to \cite{MM} for the possible occurrence of irrecoverable corners.
Finally, we remark that our weak approach might be useful to study geometric functionals depending on curvature and torsion.

\smallskip 

The uniqueness part of the proof of Theorem \ref{Tintr} 
exploits the chain rule formula in $BV$, and the expected regularity of solutions yields to an equivalent formulation of the problem \eqref{eqND-intro} with a system of integral equations. 
The existence part, instead, is based on an approximation argument, where we deal with O.D.E. systems with intitial data $\O_\e$ obtained by means of a mollification procedure applied to the distributional derivative $D\O$.  

Notice that in dimension $N=3$, if $\bG=\bI_3$ and $\O:\ol I_L\to Sk(3)$ is given by
\beq\label{Omega-intro}
\O(s):=\left(\begin{array}{ccc}
0 & \theta(s) & 0 \\ -\theta(s) & 0 & \phi(s) \\ 0  & -\phi(s) & 0  \\
 \end{array} \right)
\eeq
for some smooth functions $\t,\f \in C^1(\ol I_L)$, with $\t'(s)>0$ for every $s\in \ol I_L$, then \eqref{eqND-intro} reduces to Frenet--Serret system \eqref{FS-intro} for the given
curvature $k(s)=\t'(s)$ and torsion $\tau(s)=\f '(s)$.
More generally, according to the notation in \eqref{G-intro}, we find a unique solution to the system
\[
\left\{ \begin{array}{l}
D\gt =\gn\,D\theta \\
D\gn=-\gt\,D\theta + \gb\,D\phi \\
D\gb=-\gn\,D\phi \\
\gt(0)={\bf e}_1 \,,\quad \gn(0)={\bf e}_2\,, \quad \gb(0)= {\bf e}_3\,,
\end{array} \right. \quad s\in I_L\,.
\]
\par As a consequence, referring for the notation to Sec. \ref{Sec:notation}, we find an essentially unique rectifiable curve $\gc$ of class $C^1$ with finite total curvature and total absolute torsion given respectively by
$$
\TC(\gc)=|D\gt|(I_L)=D\t(I_L) < \i \,,\quad \TAT(\gc)=|D\gb|(I_L)=|D\f|(I_L) < \i \,.
$$
\par In order to extend the latter result to the more general case when the functions $\t$ and $\f$ in \eqref{Omega-intro} fail to be continuous, in Sec. \ref{Sec:disc} we present a case study, where we assume
$$ \t(s)=\left\{\ba{ll} s & \textrm{if } s<0  \\ s+d & \textrm{if } s>0\,,\ea\right.
\qquad
\f(s)=\left\{\ba{ll} 0 & \textrm{if } s<0  \\ \ttt & \textrm{if } s>0\,,\ea\right. 
$$
with $s\in(-1,1)$, and the jumps $d\geq 0$ and $\ttt\in\gR$ are such that $0<\sqrt{d^2+\ttt^2}<\pi$. As we shall see below, the latter bound allows us to recover the invertibility property of the solution $G(s)$ at the Jump point $s=0$.

Roughly speaking, denoting by $\gt^\pm$, $\gn^\pm$, and $\gb^\pm$ the left and right limits at $s=0$ of the tangent, normal, and binormal of the solution curve $\gc$ we obtain, it turns out that the angle between $\gn^\pm$ is
equal to $\sqrt{d^2+\ttt^2} $, whereas the angles between $\gt^\pm$ and $\gb^\pm$
are respectively given by
$$ \widehat{\gt^+\vert\gt^-}=\arccos\left(\frac{\ttt^2+d^2\cos\sqrt{d^2+\ttt^2}}{d^2+\ttt^2}\right)\,, \quad\widehat{\gb^+\vert\gb^-}=\arccos\left(\frac{d^2+\ttt^2\cos\sqrt{d^2+\ttt^2}}{d^2+\ttt^2}\right)\,. $$
\par We now observe that the datum $\O$ in \eqref{Omega-intro} can be written in terms of the canonical generators of $Sk(3)$ reported in eq. \eqref{matrix-J} through the formula
$$ \O(s)=-\t(s)\,\bJ_3-\f(s)\,\bJ_1\,. $$
Therefore, in our case study, the system in the first line of \eqref{eqND-intro} makes sense at the jump point $s=0$ provided that we replace the {\em Jump component} $D^J\O$ of the derivative of the datum with the {\em atomic measure} 
\[
D^{at}\O :=- \left( \frac {d\,\sin\sqrt{d^2+\ttt^2}}{\sqrt {d^2+\ttt^2}} \,\bJ_3 + \frac {\ttt\,\sin\sqrt{d^2+\ttt^2}}{\sqrt {d^2+\ttt^2}}\,\bJ_1\right)\,\d_0\,,
\]
where $\d_0$ is the Dirac measure at $s=0$.
\par This is the ansatz we assume in Sec. \ref{Sec:gen}, where we solve the general Frenet-Serret system in case of $BV$ data $\t$ and $\f$ in \eqref{Omega-intro}. We first consider the case of a finite jump set, and we then treat the general case by means of an approximation procedure. In any case, the essentially unique curve we obtain satisfies a bound of both its total curvature and absolute torsion in terms of the total variation of the derivative $D\O$.
\medskip\par\noindent{\bf Acknowledgments }
The research of the authors was partially supported by the GNAMPA of INDAM.
\section{Background material}\label{Sec:notation}
\subsection{Functions of bounded variations}
We refer to Secs.~3.1 and 3.2 of \cite{AFP} for the following notation.
\par Let $I\sb\gR$ be a bounded open interval, and $m\in\Nat^+$.
A vector-valued summable function $u:I\to\Rm$ is said to be of {\em bounded variation}
if its distributional derivative $Du$ is a finite $\Rm$-valued Borel measure in $I$, say $Du \in \M(I,\Rm)$.
%
The {\em total variation} $|Du|(I)$ of a function $u\in\BV(I,\Rm)$
is given by
$$ |Du|(I):=\sup\Bigl\{ \int_I \vf'(s)\,u(s)\,\dd s \mid \vf\in C^\infty_c(I,\Rm)\,,\quad \Vert\vf\Vert_\infty\leq 1 \Bigr\} $$
and hence it does not depend on the choice of the representative in the
equivalence class of the functions that agree $\calL^1$-a.e. in $I$ with $u$, where
$\calL^1$ is the Lebesgue measure in $\gR$.
Since moreover each component of $u$ is the difference of two bounded monotone functions, it turns out that
both the left and right limits $u(s^\pm):=\lim_{t\to s^\pm}u(t)$ exist for every $s\in I$, and the Jump set
$$ S(u):=\inf\{s\in I \mid u(s^-)\neq u(s^+) \} $$
is a countable set. If $I=(a,b)$, we denote $u(a)=u(a^+)$ and $u(b)=u(b^-)$.
%
%
Moreover, for $s\in I\sm S(u)$, we denote by $u(s)$ the common value $u(s^+)=u(s^-)$.
Good representatives are obtained by the left- and right-continuous functions $u_\pm(s):=u(s^\pm)$, for $s\in I$,
whereas the {\em precise} representative is obtained by choosing $u(s)=(u(s^+)+u(s^-))/2$ for every $s \in S(u)$.
\par Any function $u\in\BV(I,\Rm)$ is essentially bounded and it is differentiable $\calL^1$-a.e. on $I$, with derivative $u'$ in $L^1(I,\Rm)$.
Moreover, the distributional derivative $Du$ decomposes into the so called {\em absolutely continuous}, {\em Jump}, and {\em Cantor} parts:
$$ Du=D^au+D^Ju+D^Cu\,,\quad |Du|(I)=|D^au|(I)+|D^Ju|(I)+|D^Cu|(I)\,.$$
More precisely, one splits $Du=D^au+D^su$ into the absolutely continuous and singular parts w.r.t. the Lebesgue measure $\calL^1$.
For any $s\in S(u)$, we denote by $[u](s):=u(s^+)-u(s^-)$ the Jump of $u$ at $s$, and by $\delta_s$ the unit Dirac mass at $s$.
One has:
$$D^a u=u'\,\calL^1\,,\qquad D^Ju=\sum_{s\in S(u)}[u](s)\,\delta_s\,,\qquad D^Cu=D^su\pri(I\sm S(u))\,, $$
so that $D^Cu(B)=0$ for every countable set $B\sb I$.
We denote by $\bar Du$ the diffuse component of the derivative, i.e. $\bar Du=D^au+D^Cu$.
We thus say that $u$ is continuous if $D^Ju=0$ or, equivalently, if
$Du=\bar Du$.
\par Finally, we recall that if $u,v\in\BV(I):=\BV(I,\gR)$, the product $uv\in\BV(I)$, and the chain rule formula (cf. \cite[Sec.~3.10]{AFP}) yields:
\beq\label{chain} D^a (uv)=(u'v+uv')\,\calL^1\,,\quad D^J(uv)=\sum_{s\in S(u)\cup S(v)}[uv](s)\,\delta_s\,,\quad D^C(uv)=u D^Cv+v D^Cu \eeq
where we have set
$$ [uv](s)= u(s^+)v(s^+)-u(s^-)v(s^-)\,,
\quad s\in S(u)\cup S(v)
$$
and we can choose any good representative
of $u$ and $v$ in the third equality.
\smallskip\par A sequence $\{u_h\}\sb\BV(I,\Rm)$ is said converges to $u\in \BV(I,\Rm)$ {\em weakly-$^\ast$ in $\BV$} if $u_h$ converges to $u$ strongly in $L^1(I,\Rm)$ and
$\sup_h|Du_h|(I)<\i$. In this case, we have:
%
$$|Du|(I)\leq\liminf_{h\to\i}|Du_h|(I)\,.$$
If in addition $|Du_h|(I)\to|Du|(I)$, the sequence
$\{u_h\}$ is said to {\em strictly converge} to $u$.
The {\em weak-$^\ast$ compactness} theorem yields that if $\{u_h\}\sb\BV(I,\Rm)$ converges $\calL^1$-a.e. on $I$ to a function $u$,
and if $\sup_h|Du_h|(I)<\i$, then $u\in\BV(I,\Rm)$ and a subsequence of $\{u_h\}$ weakly-$^\ast$ converges to $u$.
%
%
%
\subsection{Curves}
In this section, we let $\gc$ be a rectifiable curve of $\RN$ with length $\calL(\gc)=L$, where $N\geq 2$.
Setting
$$I_L:=(0,L)\,,\quad \ol I_L:=[0,L]\,, \quad L>0\,,$$
the arc-length parameterization $\gc:\bar I_L\to \gR^N$ is Lipschitz-continuous, hence it is differentiable $\calL^1$-a.e. on $I_L$, by Rademacher's theorem. The tantrix $\gt$ is given by
$\gt(s):=\gc'(s)$ for $\calL^1$-a.e. $s \in I_L$, and hence $\gt\in L^1(I_L,\Sph^{N-1})$, where $\Sph^{N-1}:=\{y\in\RN\,:\,|y|=1\}$ is the Gauss hypersphere.

The {\em rotation} $K(\gp)$ of a polygonal curve $\gp$ is the sum of the turning angles at its interior vertices. Following Milnor \cite{Mi}, the {\em total curvature} $\TC(\gc)$ of $\gc$ is defined by taking the supremum of the rotation computed among inscribed polygonal curves.
Then, $\TC(\gc)<\i$ if and only if $\gt \in BV(I_L,\Sph^{N-1})$, compare e.g. \cite{Su_curv}.
In that case, moreover, one has
\[
\TC(\gc)=\textrm{Var}_{\Sph^{N-1}}(\gt)<\infty
\]
where
\[
\Var_{\Sph^{N-1}}(\gt):=\int_{I_L}|\gt'(s)|\,\dd s+|D^C\gt|(I_L)+\sum_{s\in S(\gt)}\dist_{\Sph^{N-1}}\left(\gt(s^-),\gt(s^+)\right)\,,
\]
$\dist_{\Sph^{N-1}}$ denoting the geodesic distance in $\Sph^{N-1}$.
%
Notice that
\[
\frac 2\pi\, \Var_{\Sph^{N-1}}(\gt) \leq | D\gt|(I_L)\leq \Var_{\Sph^{N-1}}(\gt)\,,
\]
where equation $| D\gt|(I_L)=\Var_{\Sph^{N-1}}(\gt)$ holds if and only if $\gt$ is continuous, i.e. when
$|D^J\gt|(I_L)=0$.

Finally, we recall that the Fr\'echet distance $d(\gc_1,\gc_2)$ between two rectifiable curves is the infimum, over all stricly monotonic
reparametrizations, of the maximum pointwise distance, see \cite[Ch. 1]{Su_curv}.
\smallskip

We now give some details of the previous facts in the physical dimension $N=3$, where we also recall from \cite{MS-AMPA} the notion of {\em total absolute torsion}.

Let $\gp$ be a polygonal curve in $\gR^3$ with consecutive vertices $v_i$, $i=0,\ldots,n$, where $n\geq 3$ and $\gp$ is not closed, i.e., $v_0\neq v_n$.
Without loss of generality, we assume that every oriented segment $\s_i:=[v_{i-1},v_{i}]$ has positive length $\calL(\s_i):=\Vert v_{i}-v_{i-1}\Vert $, for $i=1,\ldots,n$, and that two consecutive segments are never aligned, i.e., the vector product
$\s_{i}\tim \s_{i+1}\neq 0_{\gR^3}$ for each $i=1,\ldots,n-1$.
We also denote $\mesh \gp:=\sup\{\calL(\s_i)\mid i=1,\ldots,n \}$.

In the classical approach by \cite{AR,Mi}, one considers the {\em tangent indicatrix} of $\gp$, i.e., the polygonal $\ttt_\gp$ in the Gauss sphere $\Sph^2$ obtained by letting $t_i:=\s_i/\calL(\s_i)\in\Sph^2$,
for $i=1,\ldots,n$, and connecting with oriented geodesic arcs $\g_i$ the consecutive points $t_{i}$ and $t_{i+1}$, for $i=1,\ldots,n-1$.
Therefore, one has $\calL(\g_i)=d_{\Sph^2}(t_i,t_{i+1})$, and the rotation of $\gp$ is equal to the length of $\ttt_P$, i.e.
$$ K(\gp)=\sum_{i=1}^{n-1}\calL(\g_i)= \calL(\ttt_\gP)\,. 
$$
For a curve $\gc$ in $\gR^3$ as above, we thus have
$$\ba{rl} \calL(\gc):= &\sup\{\calL(\gp)\mid \gp\ll \gc\} \\ \TC(\gc):= &\sup\{ K(\gp)\mid \gp\ll \gc\}\,, \ea $$
where we write $\gp \ll \gc$ if $\gp$ is inscribed in $\gc$. The previous definitions work since for a couple of
polygonal curves $\gp$, $\gp'$ satisfying $\gp\ll \gp'$, we have $\calL(\gp)\leq \calL(\gp')$ and $K(\gp)\leq K(\gp')$.

In \cite{MS-AMPA}, we defined the total absolute torsion of $\gp$ by
$$\tilde\TAT(\gp):=\sum_{i=2}^{n-1}\wid\t_i $$
where $\wid\t_i\in[0,\p/2]$ is the shortest angle in $\Sph^2$ between the unoriented geodesic arcs $\g_{i-1}$ and $\g_i$ meeting at the edge $t_i$ of $\ttt_\gp$.
In the above notation, we are computing the sum of the distance between the consecutive discrete osculating planes (identified by their normals $\s_i\tim\s_{i+1}$). If $\gp$ is closed, a similar notation follows.
However, it may happen that for $\gp\ll \gp'$ the monotonicity inequality 
$\tilde\TAT(\gp)\leq\tilde\TAT(\gp')$ is violated.

Therefore, following the approach due to Alexandrov--Reshetnyak \cite{AR}, in \cite{MS-AMPA} we defined the total absolute torsion of $\gc$ by letting
$$
\TAT(\gc):=\lim_{\e\to 0^+}\sup\{ \tilde\TAT(\gp)\mid \gp\ll \gc \,,\,\,\m_\gc(\gp)<\e\}\,,
$$
where the {\em modulus} $\m_c(\gp)$ is given by the maximum of the diameter of the arcs of $\gc$ determined by the consecutive vertices in $\gp$.
For polygonal curves $\gp$, we have $\TC(\gp)=K(\gp)$ and $\TAT(\gp)=\tilde\TAT(\gp)$. Moreover, the following property holds:
\bp\label{PTAT} Let $\gc$ be a rectifiable curve in $\gR^3$ with both finite total curvature $\TC(\gc)$ and total absolute torsion $\TAT(\gc)$.
Then for any sequence $\{\gp_k\}$ of inscribed polygonal curves such that
$\m_\gc(\gp_k)\to 0$ as $k\to \i$, one has
$$ \lim_{k\to\i}\calL(\gp_k)=\calL(\gc)\,,\quad
\lim_{k\to\i}\TC(\gp_k)=\TC(\gc)\,,\quad \lim_{k\to\i}\TAT(\gp_k)=\TAT(\gc)\,.
$$\ep

If the curve $\gc$ is smooth, and the derivative vector $\gt'$ is non-zero everywhere on $I_L$, the scalar curvature $k:I_L\to(0,\infty)$
and the normal $\gn:I_L\to\Sph^{2}$ are identified by equation $\gt'(s)=k(s)\,\gn(s)$, the unit binormal is given by $\gb(s)=\gt(s)\tim\gn(s)$,
and the torsion $\tau:I_L\to \gr$ is computed in terms of the classical Frenet--Serret system, see eq. \eqref{F-S} below.
In addition, we have
$$
\TC(\gc)=\int_{I_L}k(s)\,\dd s\,,\quad \TAT(\gc)=\int_{I_L}|\tau(s)|\,\dd s\,.
$$

Conversely, given two smooth and bounded functions $k(s)$ and $\tau(s)$ on $I_L$, with $k(s)>0$ for every $s\in I_L$,
the fundamental theorem of curves states that there is a unique (up to rigid motions) smooth curve $\gc$ with length $L$, 
curvature $k$ and torsion $\tau$.
In fact, denoting by $({\bf e}_1,{\bf e}_2,{\bf e}_3)$ the (column) vectors of the canonical basis in $\gR^3$,
it turns out that the linear system of O.D.E.'s
\begin{equation}
\label{F-S}
\left\{ \begin{array}{l}
\gt'(s) =k(s)\,\gn(s) \\
\gn'(s)=-k(s)\,\gt(s)+\tau(s)\,\gb(s) \\
\gb'(s)=-\tau(s)\,\gn(s) \\
\gt(0)={\bf e}_1 \,,\quad \gn(0)={\bf e}_2 \,,\quad \gb(0)={\bf e}_3
\end{array} \right.
\end{equation}
has a unique solution $(\gt,\gn,\gb)$ on $I_L$, and it suffices to take $\gc(s)=0_{\gR^3}+\int_0^s \gt(\s)\,\dd\s$.
\subsection{Matrices}
\noindent
For a given integer $N\geq 2$, we denote by $\gR^{N \tim N}$ the space of square matrices of order $N$, that is equipped with the Frobenius norm 
$|A|:=\tr (A^\top A)^{1/2}$, where $\tr$ and $\top$ denote the trace and the transposition operator, respectively.
We denote by $\Lin^+_N$ the subgroup of matrices with positive determinant, and by $SO(N)$ the special orthogonal group
\[
SO(N):=\{ S\in \Lin^+_N \mid A^\top=A^{-1}\}\,,
\]
so that $\det A=1$ for every $A\in SO(N)$. Moreover, $\bI=\bI_N$ is the identity matrix and $Sk(N)$ the Lie algebra of skew-symmetric matrices, i.e.
\[  Sk(N) :=\left\{ B\in \gR^{N\tim N} \mid B^\top=-B\right\}\,.
\]
For any $B \in Sk(N)$, the square matrix $B^2$ is symmetric, whereas the exponential $A:=\exp B$ is an orthogonal matrix in $SO(N)$.
Conversely, every $A\in SO(N)$ can be written as $A=\exp B$ for some $B \in Sk(N)$.
%

In the physical dimension $N=3$, a non-trivial matrix $B\in Sk(3)$ is given by $B=\a\,\bJ$ for some positive number $\a>0$ and some
$\bJ\in Sk(3)$ such that $\bJ^3=-\bJ$ and $\bJ^4=\bI_3$, as e.g. in definition \eqref{matrix-J} below.
Rodrigues formula for the matrix $A=\exp B\in SO(3)$ gives
\beq\label{Rodrigues}
A=\exp(\a\,\bJ)=\bI_3+\sin\a\,\bJ+(1-\cos\a)\,\bJ^2
\eeq
and $A$ is a rotation with angle $\a$ and axis oriented by the axial unit vector $\gk\in\Sph^2$ of $\bJ$,
so that
\[ \gk\tim\gv=\bJ\gv\,,\quad \fa\,\gv\in\gR^3\simeq\gR^{3\tim 1}.
\]
In particular, when $\gk={\bf e}_\ell$, where $\ell=1,2,3$,
the previous vector product formula is satisfied with $\bJ_\ell$ given respectively by
\beq\label{matrix-J}
\bJ_1 := \left(
     \begin{array}{ccc}
       0 & 0 & 0 \\ 0 & 0 & -1 \\ 0 & 1 & 0 \\      \end{array}   \right)\,,
       \quad
\bJ_2 := \left(
     \begin{array}{ccc}
       0 & 0 & 1 \\ 0 & 0 & 0 \\ -1 & 0 & 0 \\      \end{array}   \right)\,,
       \quad
\bJ_3 := \left(
     \begin{array}{ccc}
       0 & -1 & 0 \\ 1 & 0 & 0 \\ 0 & 0 & 0 \\      \end{array}   \right)\,,
\eeq
where the triplet $(\bJ_1,\bJ_2,\bJ_3)$ is a generator of the Lie algebra $Sk(3)$.

\smallskip

For our purposes, in the first appendix we report a proof of the following result that goes back to Cayley.
%
\bp\label{P-Cayley} Let $A\in SO(N)$ be such that $A+\bI_N$ is invertible.
Then 
$$B:=(A+\bI_N)^{-1}(A-\bI_N) \in Sk(N)\ . $$
\ep
In low dimension $N=2$, if $A\in SO(2)$, condition $\det(A+\bI_2)=0$ is violated if and only if $A$ represents a planar rotation with angle $\pi$, see Example \ref{E2D} below.

In dimension $N=3$, if $A\in SO(3)$ is given by Rodrigues formula \eqref{Rodrigues}, we have $\det(A+\bI_3)\neq 0$ if and only if $\cos\a\neq -1$. For $|\a|<\pi$, we have
\[
B:=(A+\bI_3)^{-1}(A-\bI_3)=\left( \bI_3+\sin\a\,\bJ+(1-\cos\a)\,\bJ^2
\right)^{-1}\left( \sin\a\,\bJ+(1-\cos\a)\,\bJ^2
\right)\,.
\]
Choosing e.g. $\bJ=\bJ_2$, see eq. \eqref{matrix-J},
we get
\[ A:=\exp(\a\,\bJ_2)=\left(
     \begin{array}{ccc}
       \cos\a & 0 & \sin\a \\ 0 & 1 & 0 \\ -\sin\a & 0 & \cos\a \\      \end{array}   \right)
\]
and hence
\[
 A-\bI_3=\left(
     \begin{array}{ccc}
       \cos\a-1 & 0 & \sin\a \\ 0 & 0 & 0 \\ -\sin\a & 0 & \cos\a-1 \\      \end{array}   \right)\,,
\quad
A+\bI_3=\left(
     \begin{array}{ccc}
       \cos\a+1 & 0 & \sin\a \\ 0 & 2 & 0 \\ -\sin\a & 0 & \cos\a+1 \\      \end{array}   \right)
\]
so that $\det(A+\bI_3)=4(\cos\a+1)\neq 0$ for $|\a|<\pi$, and in that case
\[
(A+\bI_3)^{-1}=\frac 12\left(
     \begin{array}{ccc}
       1 & 0 & -\l \\ 0 & 1 & 0 \\ \l & 0 & 1 \\      \end{array}   \right)\,,\quad
\l:=\frac{\sin\a}{\cos\a+1}=\tan(\a/2)\,.
\]
We thus obtain
\[
B:=(A+\bI_3)^{-1}(A-\bI_3)=\tan(\a/2)\,\bJ_2\,,
\quad |\a|<\pi\,.
\]
\subsection{Matrix valued BV functions}
Let $G:I_L\to \gR^{N\tim N}$ be a matrix-valued function in $BV(I_L,\gR^{N\tim N})$, so that the distributional derivative of $G$ is a matrix valued Borel measure $DG \in \M(I_L,\gR^{N\tim N})$ such that $|DG|(I_L)<\infty$.
The precise representative of $G$ is given by
\beq\label{G-pr}
G(s)=\frac 12\, \left( G(s^+)+G(s^-) \right)\,,\quad s \in I_L
\eeq
where $G(s^\pm)$ denote the right and left limits, and we let $G(0)=G(0^+)$, $G(L)=G(L^-)$.
Moreover, the Jump set of $G$ is the countable set
\[
S(G)=\{s \in I_L \mid G(s^-)\neq G(s^+)\}\,.
\]
For $X=\Lin^+_N$, $SO(N)$, or $Sk(N)$, we say that $G\in BV(I_L,X)$ if in addition $G(s)\in X$ for $\calL^1$-a.e. $s\in I_L$.
\smallskip

If e.g. $N=3$ and $G(s)=\exp(\t(s)\bJ_2)$ for some smooth angle function $\t:\ol I_L\to\gR$, where $\bJ_2\in Sk(3)$ is given by \eqref{matrix-J}, we readily obtain
\[
G'(s)= \left(\begin{array}{ccc}
\sin\t(s) & 0 & \cos\t(s) \\
0 & 0 & 0 \\
-\cos\t(s) & 0 & -\sin\t(s) \\
\end{array}\right)\,\theta'(s)\quad \fa\, s \in \ol I_L\,.
\]

Let $G\in BV(I_L,SO(N))$, and assume that $S(G)=\emptyset$, so that $G(s)$ is continuous on $I_L$. Since $G(s)^\top G(s)=\bI_N$ for every $s$, by the chain rule formula we can write
\[
0=D(G^\top G)=DG^\top G+G^\top DG\,.
\]
Therefore, letting $B:=DG^\top G$, we have $B\in \M(I_L,\gR^{N\tim N})$, and the absolutely continuous component of $B$ has a density function that belongs to $L^1(I_L,Sk(N))$. We also have $-B=B^\top=G^\top DG$, so that since $G(s)$ is invertible for every $s$ we obtain
\[
DG=-(G^\top)^{-1}B=-G\,B\,.
\]

If $D^JG\neq 0$, instead, for every $s\in S(G)$ the precise representative satisfies equation
\[ G(s)=\frac 12\,G(s^-)\left( A(s)+\bI_N\right)
\]
where
\beq\label{As}
A(s):=G(s^-)^\top G(s^+)\in SO(N)\,.
\eeq
Therefore, $A(s)=\bI_N$ if $s\in I_L\sm S(G)$, and $G(s)$ is invertible if and only if
\beq\label{detAs}
\det\left( A(s)+\bI_N\right)\neq 0\,.
\eeq
Notice that
eq. \eqref{detAs} is violated if and only if $-1$ is an eigenvalue of $A(s)$.
If eq. \eqref{detAs} holds for every $s\in S(G)$,
the atomic part $B^{at}$ of the measure $B:=DG^\top G$ is given by
\[
B^{at}=\sum_{s\in S(G)}B^{at}(s)\,\d_s
\]
where for any $s\in S(G)$ we have set
\[
B^{at}(s)=-\frac 12\,\left( G(s^+)^\top+G(s^-)^\top \right)\left( G(s^+)-G(s^-) \right)=\frac 12\left( A(s)^\top-A(s)\right)\in Sk(N)\,.
\]
\bex\label{E2D}
If e.g. $N=2$ and for some $s_0\in S(G)$
\[
G(s_0^\pm)=\left(\begin{array}{cc}
       \cos\a_\pm & -\sin\a_\pm \\
       \sin\a_\pm & \cos\a_\pm \\
     \end{array} \right)\,,
\]
setting $\DD\a=\a_+-\a_-$ we have
\[ A(s_0)=\left(\begin{array}{cc}
       \cos\DD\a & -\sin\DD\a \\
       \sin\DD\a & \cos\DD\a \\
     \end{array} \right)\,,
\quad
B^{at}(s_0)=\left(\begin{array}{cc}
       0 & \sin\DD\a \\
       -\sin\DD\a & 0 \\
     \end{array} \right)\,.
     \]
Therefore, eq. \eqref{detAs} is violated if and only if $\cos\DD\a=-1$, i.e. $A(s_0)=-\bI_2$.
\eex
%
%
\section{The case of continuous data}\label{Sec:cont}
In this paper, we wish to extend the fundamental theorem of curves in the Euclidean space to a wider range of data.
To this purpose, we first consider a general problem in the framework of matrix valued functions of bounded variation, in any dimension $N\geq 2$.
The main issue is extending a standard result on linear O.D.E. systems to the case 
when the datum is a matrix valued measure.
\subsection{A general problem}
Let $\O : I_L\to Sk(N)$ be a continuous function in $BV(I_L,Sk(N))$, 
where $I_L=(0,L)$ and $L>0$, and let $\bG \in SO(N)$ be a given initial condition. 
We look for
solutions $G$ in $BV(I_L,\gR^{N\tim N})$ to the system
\begin{equation}
\label{eqND}
\left\{ \begin{array}{l}
DG= - (G^\top)^{-1}D\O \\
G(0)=\bG\,.
\end{array} \right.
\end{equation}
\par We remark that in the first line of the system \eqref{eqND} we have an equation between measures, and it has to be intended by taking the precise representative of a function $G$ in $BV(I_L,\gR^{N\tim N})$.
Therefore, we are looking for solutions such that $G(s)$ is an invertible matrix in $\gR^{N\tim N}$ for every $s\in I_L$.
In this framework, we shall see that solutions to \eqref{eqND} are continuous functions taking values in $SO(N)$.

More precisely, denoting by $H^i_j(s)$ the generic entry of a matrix valued function $H(s)$, the given equation is satisfied provided that for every $i,j=1,\ldots,N$
and for every test function $\vf\in C^\i_c(I_L)$ we have:
\beq\label{eq-weak}
\int_{I_L}G^i_j(s)\,\vf'(s)\,\dd s=\int_{I_L}\vf(s)\,\left(
(G(s)^\top)^{-1} \,\dd D\O(s)\right)^i_j\,.
\eeq

Notice that in dimension $N=3$, if $\bG=\bI_3$ and $\O:\ol I_L\to Sk(3)$ is given by
\begin{equation}
\label{Omega3}
\O(s):=\left(\begin{array}{ccc}
0 & \theta(s) & 0 \\ -\theta(s) & 0 & \phi(s) \\ 0  & -\phi(s) & 0  \\
 \end{array} \right)
\end{equation}
for some smooth functions $\t,\f \in C^1(\ol I_L)$, with $\t'(s)>0$ for every $s\in \ol I_L$, then \eqref{eqND} reduces to the classical Frenet--Serret system \eqref{F-S} for the given
curvature $k(s)=\t'(s)$ and torsion $\tau(s)=\f '(s)$, that is well-known to have a unique solution $G(s)\in C^1(\ol I_L,SO(3))$.
\smallskip

The following general result ensures the uniqueness and the expected regularity of solutions to the system \eqref{eqND}, provided that it admits solutions.
\bp\label{T-uniqueness}
For a given continuous map $\O\in BV(I_L,Sk(N))$ and $\bG\in SO(N)$, assume that the system
\eqref{eqND} admits a solution in $BV(I_L,\gR^{N\tim N})$. Then, the solution is unique, and it is a continuous map in the class $BV(I_L,SO(N))$.
In particular, \eqref{eqND} is equivalent to the system
\begin{equation}
\label{eqND-cont}
\left\{ \begin{array}{l}
DG= - G\,D\O \\
G(0)=\bG\,.
\end{array} \right.
\end{equation}
\ep
\bpf If $G\in BV(I_L,\gR^{N\tim N})$ solves the system \eqref{eqND}, by the equivalent equation $G^\top DG =-D\O$ between measures, where $D^J\O=0$, we obtain that $D^JG=0$, and hence that $G$ is continuous. Furthermore, we have
\[ DG^\top G=(G^\top DG )^\top=-(D\O)^\top=D\O\,,
\]
whence by the chain rule formula
\[
D(G^\top G)=DG^\top G+G^\top DG=D\O-D\O=0
\]
and hence $G^\top G$ is constant on the open interval $I_L$.
Since $\bG\in SO(N)$, by the initial condition we have $G(0)^\top G(0)=\bG^\top\bG=\bI_N$. Therefore, we infer that $G(s)^\top G(s)=\bI_N$ and hence $G(s)\in SO(N)$ for every $s\in I_L$, so that $(G^\top)^{-1}=G$. Assume now that $\wid G\in BV(I_L,\gR^{N\tim N})$ is another solution to \eqref{eqND}, and let $H:=G-\wid G$. Then, $H\in BV(I_L,\gR^{N\tim N})$ is continuous, and by linearity
$$ DH=D(G-\wid G)=-G\,D\O+\wid G\,D\O=-H\,D\O \Lra H^\top DH=-D\O\,. $$
Therefore, recalling that the measure $D\O$ is skew-symmetric, we have
$$ DH^\top H=(H^\top DH)^\top=(-D\O)^\top=D\O\,,$$
and the chain-rule formula gives
\[
D(H^\top H)=H^\top DH+DH^\top H= -D\O+D\O=0
\]
and hence the continuous function $H^\top H$ is constant on $I_L$. Since $H(0)=G(0)-\wid G(0)=0$, we have $(H^\top H)(0)=0$, so that $H^\top H=0$ on $I_L$.
Therefore, by equation $|H|^2=\tr\left(H^\top H\right)$ we deduce that $H\equiv 0$ on $I_L$ and hence that the solution to \eqref{eqND} is unique. 
The last assertion readily follows.
\epf
\par
As a consequence, we have:
\bp\label{P-integral-eq}
For a given continuous map $\O\in BV(I_L,Sk(N))$ and $\bG\in SO(N)$, the system \eqref{eqND} is equivalent to the integral equation
\beq\label{integral-eq}
G(s)-\bG=-\int_0^s G(\s)\,\,\dd D\O(\s)\,,\quad \fa\,s\in I_L
\eeq
among continuous functions $G\in BV(I_L,SO(N))$. 
\ep
\bpf
By Proposition \ref{T-uniqueness}, a solution $G$ to \eqref{eqND} is a continuous 
function $G\in BV(I_L,SO(N))$ solving eq. \eqref{eqND-cont}.
Therefore, for every $s\in I_L$ we have
$$ G(s)-G(0)=\int_0^s \dd\,DG(\s) $$
and hence $G(s)$
satisfies the integral equation \eqref{integral-eq}. Conversely, we show that if 
$G\in BV(I_L,SO(N))$ is a continuous solution to \eqref{integral-eq}, then eq.
\eqref{eq-weak} holds for every $i,j=1,\ldots,N$ and $\vf \in C^\i_c(I_L)$.
Now, eq. \eqref{integral-eq} in components reads
$$
G(s)^i_j-\bG^i_j =-\int_0^s
G(\s)^i_k\,\dd D\O(\s)^k_j \,,\quad \fa\,s\in I_L\,,
$$
where the summation on the repeated index $k=1,\ldots,N$ is implicitly assumed. We thus have
$$
\ds \int_{I_L}G(s)^i_j\,\vf'(s)\,\dd s =
\int_{I_L}\vf'(s)\,\bG^i_j\,\dd s -\int_{I_L}\vf'(s)\left( \int_0^s G(\s)^i_k\,\dd D\O(\s)^k_j\right)\,\dd s \,,
$$
where the first addendum in the right-hand side is equal to zero.
Since moreover both the Borel measures $\calL^1$ and $D\O$ are $\sigma$-finite, by Fubini-Tonelli's theorem we have
$$
\begin{array}{rl}
\ds -\int_{I_L}\vf'(s)\Bigl( \int_0^s G(\s)^i_k\,\dd D\O(\s)^k_j\Bigr)\,\dd s = &
\ds -\int_{I_L}\Bigl( \int_\sigma^L \vf'(s)\,\dd s\Bigr) G(\s)^i_k\,\dd D\O(\s)^k_j \\
= & \ds \int_{I_L} \vf(\s)\,G(\s)^i_k\,\dd D\O(\s)^k_j  \\
= & \ds \int_{I_L}\vf(s)\,\left(G(s)\,\dd D\O(s)\right)^i_j
\end{array}
$$
and hence \eqref{eq-weak} holds, since $G(s)=(G(s)^\top)^{-1}$ for every $s\in I_L$.
\epf

The following existence theorem is the first main result of this paper.
\bt\label{T-existence}
For any given continuous map $\O\in BV(I_L,Sk(N))$ and $\bG\in SO(N)$, the system
\eqref{eqND} admits a solution in $BV(I_L,\gR^{N\tim N})$.
\et
\bpf
It is divided in four steps.
\smallskip

\noindent{\em Step 1: approximation.}
We extend $\O(s)$ to the neighborhood $I=(-\e_0,L+\e_0)$ of $\ol I_L$, with $\e_0>0$ small, by letting $\O(s)=\O(0)$ for $s<0$, and $\O(s)=\O(L)$ for $s>L$.
Let $\r\in C^\i_c(\gR)$ be a smooth symmetric mollifier, i.e. $\r$ is non negative, even, with support contained in $(-1,1)$ and integral $\int_{-1}^1 \r(s)\,\dd s=1$. For $0<\e<\e_0$, we let $\r_\e(s):=\e^{-1}\r(\e^{-1}s)$ and we define $B_\e:I\to Sk(N)$ by means of convolutions
of the derivative measure $D\O$, i.e.
$$ B_\e(s):=\int_{I}\r_\e(s-\s)\,\dd\,D\O(\s)\,,
\quad s\in I\,. $$
Therefore, for every $i,j=1,\ldots, N$ we have
\beq\label{tv-conv}
\lim_{\e\to 0^+} 
\int_{I}|B_\e(s)^i_j|\,\dd s=|D\O^i_j|(I_L)\,.
\eeq
%
%
Since moreover the function $B_\e(s)$ is smooth, the problem
\beq\label{eq3D-e}
\left\{ \begin{array}{l}
H'(s)= - H(s)\,B_\e(s) \\
H(0)=\bG
\end{array} \right.
\quad s\in I\,,
\eeq
has a unique solution $G_\e\in C^\i(I,SO(N))$,
that equivalently solves the integral equation
\[
H(s)=\bG-\int_0^s H(\s)B_\e(\s)\,\dd\s\,,\quad s\in I\,.
\]
Since $G_\e(s)$ is an orthogonal matrix, for every $s\in I$ we have
\[
|G_\e'(s)|= |G_\e(s)\,B_\e(s)|=|B_\e(s)|
\]
and hence by the properties of convolution products we can estimate
\[
\int_{I}|G_\e'(s)|\,\dd s\leq |D\O|(I_L)<\i\,.
\]
As a consequence, for every sequence $\{\e_h\}\searrow 0$ we can find a (not relabeled) subsequence and a function $G\in BV(I,\gR^{N\times N})$ such that
$\{G_{\e_h}\}$ weakly-* converges in $BV$ to $G$.
\smallskip

\noindent{\em Step 2: properties of the limit function.}
We show that $G$ is a continuous map in $BV(I,SO(N))$ such that $G(0)=\bG$.

In fact, by the pointwise $\calL^1$-a.e. convergence of a not relabeled subsequence of $\{G_{\e_h}\}$ to $G$, we infer that $G(s)\in SO(N)$ for $\calL^1$-a.e. $s\in I$. Moreover, if $-\e_0<a<b<L+\e_0$ are chosen in such a way that $|D^JG|(\{a\})=|D^JG|(\{b\})=0$, by lower semicontinuity we have
$$ |DG|((a,b))\leq\liminf_{h\to\i}\int_a^b|G_{\e_h}'(s)|\,\dd s\leq |D\O| ((a,b))\,.
$$
Since for every $s_0\in I$ we can find a sequence $r_k\searrow 0$ such that
$|D^JG|(\{s_0\pm r_k\})=0$ for every $k$, recalling that the $Sk(N)$-valued measure $D\O$ has no atomic part we obtain
$$
\limsup_{k\to\i}|DG|((s_0-r_k,s_0+r_k))\leq \limsup_{k\to\i} |D\O|((s_0-r_k,s_0+r_k))=0\,.
$$
Therefore, $D^JG=0$ and hence the function $G(s)$ has a continuous representative.
In particular, the sequence $\{G_{\e_h}(s)\}$ converges to $G(s)$ for every $s\in I$. Since $G_{\e_h}(0)=\bG$ for every $h$, we obtain $G(0)=\bG$.
\smallskip

\noindent{\em Step 3: existence.}
By Proposition \ref{P-integral-eq}, it suffices to show that the function $G$ this way obtained solves the integral equation \eqref{integral-eq}.
%
Recalling that the function $G_{\e_h}$ satisfies equation
$$
G_{\e_h}(s)-\bG=-\int_0^s G_{\e_h}(\s)\,B_{\e_h}(\s)\, d\s\,,\quad \fa\,s\in I_L\,,
$$
and that $G_{\e_h}(s)\to G(s)$ for every $s \in I_L$, it suffices to show that for every $i,j,k=1,\ldots,N$ we have
$$ \lim_{h\to\i}\int_0^{s_0} G_{\e_h}(s)^i_k\,B_{\e_h}(s)^k_j\,\dd s=
\int_0^{s_0} G(s)^i_k\,\dd\, D\O(s)^k_j \quad \fa\,s_0\in I_L\,. $$
\par The latter formula holds true by applying the following result on measure theory, the proof of which is given in Step 4.
\bl\label{L-cont} Let $A\sb \gR$ be a bounded interval and let $\{\m_h\}$ be a sequence of signed Borel measures in $A$ weakly-$^*$ converging to a signed Borel measure $\m$ in $A$ and such that $|\m_h|(A)\to|\m|(A)$ as $h\to\i$. Let $\{f_h\}$ be a sequence of continuous functions on $A$ converging pointwise $\m$-a.e. in $A$ to some continuous function $f$, and such that $\sup_h\Vert f_h\Vert_\i\leq C<\i$. Then
$$ \lim_{h\to\i}\int_A f_h\,\dd\m_h=\int_A f\,\dd\m\,.
$$
\el
\par
In fact, letting $f_h(s)=G_{\e_h}(s)^i_k$, $f(s)=G(s)^i_k$, $\m_h= (B_{\e_h})^k_j\,\calL^1$, $\m=(D\O)^k_j$, and $A=(0,s_0)$, the hypotheses of Lemma \ref{L-cont} are verified due to the pointwise convergence of $\{G_{\e_h}\}$ to $G$ and to the weak-$^*$
convergence $B_{\e_h}\,\calL^1 \wc D\O$ joined with $|B_{\e_h}\,\calL^1|(I_L)\to|D\O|(I_L)$.
\smallskip

\noindent{\em Step 4: proof of Lemma \ref{L-cont}.} Let $\e>0$. By Egorov theorem, there exists a Borel set $K_\e\sb A_\e$ such that $\m(A_\e\sm K_\e)<\e$ and the sequence $\{f_h\}$ converges to $f$ uniformly on $K_\e$. We have
$$ \ba{crl} \ds\Bigl\vert \int_A f_h\,\dd\m_h - \int_A f\,\dd\m \Bigr\vert 
& \leq  &
\ds \Bigl\vert \int_{K_\e} f_h\,\dd\m_h - \int_{K_\e} f\,\dd\m \Bigr\vert
+
\int_{A \sm K_\e} |f| \,\dd|\m|
+
\int_{A \sm K_\e} | f_h | \,\dd|\m_h| \\
& =: & I_1+I_2+I_3\,. 
\ea $$
By the uniform convergence in $K_\e$ and by the weak-$^*$ convergence $\m_h\wc \m$ as measures, we infer that $I_1<\e$ for $h$ large enough.
Since moreover by the pointwise convergence $\Vert f \Vert_\i\leq C$, we have that
$$
I_2\leq
C \m(A \sm K_\e)\leq C\,\e\,. $$
Finally, since by the total variation convergence
$|\m_h|(A \sm K_\e)\leq 2\,
|\m| (A \sm K_\e)$ for $h$ large, we similarly estimate
$$
I_3 \leq 
C |\m_h|(A \sm K_\e)\leq 2\,C\,\e $$
and the assertion readily follows. 
\epf
\subsection{The weak Frenet-Serret system with continuous data}
We  now consider the system \eqref{eqND} in dimension $N=3$, where we assume that both the angle functions $\t$ and $\f$
of the $Sk(3)$-valued $BV$ function in \eqref{Omega3} are continuous.

We recall that in this case \eqref{eqND} is equivalent to the system \eqref{eqND-cont}.
Therefore, letting
\beq\label{tnb}
\gt(s)=G(s){\bf e}_1\,,\quad \gn(s)=G(s){\bf e}_2\,,\quad \gb(s)=G(s){\bf e}_3 \,,
\eeq
our problem \eqref{eqND} is equivalent to
\[
\left\{ \begin{array}{l}
D\gt =\gn\,D\theta \\
D\gn=-\gt\,D\theta + \gb\,D\phi \\
D\gb=-\gn\,D\phi \\
\gt(0)={\bf e}_1 \,,\quad \gn(0)={\bf e}_2\,, \quad \gb(0)= {\bf e}_3
\end{array} \right.
\]
so that if $\theta,\phi \in W^{1,1}(I_L)$, letting $k(s):=\theta'(s)$ and $\tau(s):=\phi'(s)$, we recover the Frenet--Serret system
\eqref{F-S} in the $\calL^1$-a.e. sense.

If $\O$ is smooth, and $\theta$ is strictly increasing, we thus obtain the classical fundamental theorem of curves.
In the second appendix we discuss to which extension one can find the explicit formula of solutions to the classical Frenet--Serret system.

As a consequence of our previous results, we readily obtain the following
\bt\label{Tmain} Let $N=3$, and assume that the datum $\O\in BV(I_L,Sk(3))$ is given by \eqref{Omega3}, where both the angle functions $\t$ and $\f$ are continuous.
Then, the system \eqref{eqND} admits a unique continuous solution $G\in BV(I_L,SO(3))$.
Moreover, according to \eqref{tnb}, the Frenet--Serret system \eqref{F-S} holds true for $\calL^1$-a.e. $s \in I_L$, where $k(s)=\t'(s)$ and $\tau(s)=\f'(s)$.
Finally, if the function $\theta$ is strictly increasing, there exists a unique (up to rigid motions) rectifiable curve $\gc:\bar I_L\to\gR^3$ of class $C^1$, parameterized in arc-length, with tantrix
$ \gc'(s)=\gt(s)$ for every $s\in I_L$, and with finite total curvature
$$ \TC(\gc)=D\t(I_L)=\theta(L)-\theta(0)\,. $$
\et
%
\bpf 
The first two assertions follow from Proposition \ref{T-uniqueness} and Theorem \ref{T-existence}.  In addition, if $\theta$ is strictly increasing, the given curve is essentially given by
$$ \gc(s):=0_{\gR^3}+\int_0^s \gt(\s)\,\dd\s\,,\quad s\in \ol I_L\,, $$
and $\gc$ is the arc-length parameterization of a regular curve of class $C^1$.
Since moreover $\gt:I_L\to\Sph^2$ is continuous and with bounded variation, we already know that
$$ \TC(\gc)=|D\gt|(I_L)<\i $$
and hence the last assertion follows since $|D\gt|(I_L)=D\t(I_L)$.
\epf

We now deal with the total absolute torsion of the solution curve:
\bt\label{T-tor}
Under the hypotheses of Theorem $\ref{Tmain}$, the curve $\gc$ has finite total absolute torsion given by
$$
\TAT(\gc)=|D\gb|(I_L)=|D\f|(I_L) < \i \,.
$$
\et
\bpf
Coming back to the proof of Theorem \ref{T-existence}, where $N=3$ and $\O(s)$ is given by \eqref{Omega3}, we have
$$ B_\e(s)=
\left(\begin{array}{ccc}
0 & k_\e(s) & 0 \\ -k_\e(s) & 0 & \tau_\e(s) \\ 0  & -\tau_\e(s) & 0  \\
 \end{array} \right)
$$
where
$$ k_\e(s):=\int_{I}\r_\e(s-\s)\,\dd\,D\t(\s)\,,\quad \tau_\e(s):=\int_{I}\r_\e(s-\s)\,\dd\,D\f(\s)\,,\quad s\in I\,. $$
For every $h\in\Nat$ we let
$$ \gc_h(s):=0_{\gR^3}+\int_0^s \gt_h(\s)\,\dd\s\,,\quad s\in \ol I_L\,, $$
where the ordered triplet $(\gt_h,\gn_h,\gb_h)$ is given by the column vectors of the approximate solution $G_{\e_h}(s)\in C^\i(\ol I_L,SO(3))$. Whence the rectifiable curve $\gc_h$ satisfies
$$ \TAT(\gc_h)=\int_{I_L}|\gb'_h(s)|\,\dd s=\int_{I_L}|\tau_{\e_h}(s)|\,\dd s $$
for every $h\in\Nat$, so that by the total variation convergence
in \eqref{tv-conv}
we infer that
\beq\label{TATch}
\lim_{h\to\i}\TAT(\gc_h)=\lim_{h\to\i}\int_{I_L}|\tau_{\e_h}(s)|\,\dd s= |D\f|(I_L)<\i\,.
\eeq

In addition, for every interval $(a,b)\sb I_L$ we have
$$ \int_a^b|\gt_h'(s)|\,ds=\int_a^b k_{\e_h}(s)\,ds\leq |D\theta|((a,b)) \quad\fa\,h\,,$$
where the measure $D\theta$ has no atomic part, so that the sequence $\{\gt_h\}$ is equibounded and equi-absolutely continuous. Therefore, by Ascoli--Arzel\`a theorem a (not relabeled) subsequence of $\{\gt_h\}$ converges to $\gt=\dot\gc$ uniformly, and hence the corresponding sequence $\{\gc_h\}$ converges to $\gc$ in the $C^1$-sense.
\par For $k\in\Nat$ large, we consider the subdivision of $\ol I_L$ in $k$ equal segments, i.e. we let $s_n^k=(L/k)\,n$, where $n=0,\ldots,k$. We denote by
$\gp_k$ the curve inscribed in $\gc$ and obtained
by consecutively connecting with straight segments the vertices $\gc(s^k_n)$. For every $h$, the polyhedral curve $\gp^{(h)}_k$ inscribed in $\gc_h$ is obtained in a similar way, so that its consecutive vertices are the points $\{\gc_h(s^k_n)\}_{n=0}^k$.

We shall denote by $\ttt_\gp$ the tantrix of a polyhedral curve $\gp$, and by $\bbb_\gp$ the corresponding polar curve. 
Notice that by the $C^1$ convergence of $\{\gc_h\}$ to $\gc$, we can find $m_0\in\Nat$ such that for $k\geq m_0$ the turning angles of the polygonals $\gp^{(h)}_k$ are small, for every $h\geq m_0$, so that we definitely obtain
\beq\label{TAT-pol}
\TAT(\gp^{(h)}_k)=\TC_{\Sph^2}(\ttt_{\gp^{(h)}_k})=\calL(\bbb_{\gp^{(h)}_k})
\eeq
for every $h,\,k\geq m_0$,
where $\TC_{\Sph^2}$ denotes the total intrinsic curvature in $\Sph^2$.

Since moreover $\m_\gc(\gp_k)\to 0$, by the $C^1$ convergence of $\{\gc_h\}$ to $\gc$ we infer that $\m_{\gc_h}(\gp^{(h)}_k)\to 0$ as $k\to \i$ for every $h\in\Nat$.
Therefore, by Proposition \ref{PTAT} we have
\beq\label{lim1}
\lim_{k\to\i}\TAT(\gp^{(h)}_k)=\TAT(\gc_h)\quad\fa\,h\,.
\eeq
\bl\label{L-lim} We also have
$\ds\lim_{h\to\i}\TAT(\gp^{(h)}_k)=\TAT(\gp_k)$ for every $k$ large.
\el
\bpf 
Every curve $\gc_h$ is parameterized in arc-length and $\{\gc_h\}$ converges to $\gc$ in the $C^1$ sense. We thus infer that for any fixed $k$ the vertices $\gc_h(s_n^k)$ of the polygonal $\gp^{(h)}_k$ converge to the vertices $\gc(s_n^k)$ of $\gp_k$, whence $\calL(\gp^{(h)}_k)\to\calL(\gp_k)$ as $h\to\i$.
Furthermore, each angle between two consecutive vertices of $\gp^{(h)}_k$ converges to 
the corresponding angle of $\gp_k$. Therefore, the corresponding polar curves satisfy $\calL(\bbb_{\gp^{(h)}_k})\to\calL(\bbb_{\gp_k})$ as $h\to\i$, and the assertion follows since $\TAT(\gp^{(h)}_k)=\calL(\bbb_{\gp^{(h)}_k})$ and $\TAT(\gp_k)=\calL(\bbb_{\gp_k})$  for every $k$ and $h$.
\epf
\par Now, using Lemma \ref{L-lim}, \eqref{lim1} and \eqref{TATch} we can estimate
$$
\ba {rl}
\ds\limsup\lm_{k\to\i} \TAT(\gp_k)= & \limsup\lm_{k\to\i}\lim\lm_{h\to\i}\TAT(\gp^{(h)}_k) \\
\leq & \limsup\lm_{h\to\i}\lim\lm_{k\to\i}\TAT(\gp^{(h)}_k)=
\limsup\lm_{h\to\i}\TAT(\gc_h)= |D\phi|(I_L) \ea $$
and hence the curve $\gc$ has finite total absolute torsion bounded by
$$ \TAT(\gc)\leq |D\phi|(I_L)<\i \,. $$
Since moreover by Proposition \ref{PTAT} we know that $\TAT(\gc)$ is equal to the limit of the sequence $\{\TAT(\gp_k)\}$, using \eqref{TATch}, \eqref{lim1}, and Lemma \ref{L-lim} we can estimate
$$
\ba{rl} |D\phi|(I_L)= &
\ds \lim_{h\to\i}\TAT(\gc_h)=
\lim_{h\to\i}\lim_{k\to\i} \TAT(\gp^{(h)}_k)\leq \limsup_{k\to\i}\limsup_{h\to\i} \TAT(\gp^{(h)}_k) \\
= &
\ds \limsup_{k\to\i} \TAT(\gp_k)=\TAT(\gc)\,.
\ea
$$
Therefore, recalling that $|D\gb|(I_L)=|D\phi|(I_L)$, the assertion readily follows.
\epf
\section{The case of discontinuous data}\label{Sec:disc}
We wish to extend the result in Theorems \ref{Tmain} and \ref{T-tor}
to the case of discontinuous data.
In order to focus the problem under scrutiny, we first briefly discuss the easier case of dimension $N=2$.
We then give an example of a solution curve with a corner point where both the derivatives of the functions $\t$ and $\f$ in \eqref{Omega3} have a non-trivial atom.
The analysis at the singular point suggests what we believe is a geometrically consistent approach to extend eq. \eqref{eqND} to the case of a $BV$ datum $\O$ with jump points.
\medskip

Returning to the general problem, we now assume that the datum $\O\in BV(I_L,Sk(N))$ fails to be continuous, i.e. $D^J\O\neq 0$. Then, the system \eqref{eqND} cannot have continuous solutions $G\in BV(I_L,\gR^{N\tim N})$. If it were the case, in fact, the function $G$ would solve the equivalent equation $-G^\top DG=D\O$, and this is not possible. Therefore, we have to look for a reasonable way to give sense to the system \eqref{eqND}.

Notice that in the easier case
when the Jump set of $\O$ is finite, the equation in the first line of \eqref{eqND}
makes sense separately on each connected component of $I_L\setminus S(\O)$, and hence it is reasonable to look for solutions $G\in BV(I_L,\Lin^+_N)$ that are continuous on $I_L\setminus S(\O)$, so that $S(G)\sb S(\O)$.
However, some care has to be taken at the discontinuity points of the datum $\O$.

\par In fact, if we consider the left continuous representative
\[
G_-(s):=G(s^-)\,,\quad s \in I_L
\]
of a function $G\in BV(I_L,SO(N))$, the system \eqref{eqND} is equivalent to
\beq
\label{sist-wrong}
\left\{ \begin{array}{l}
DG= -G_-M D\O \\
G(0)=\bG
\end{array} \right.
\eeq
where on account of \eqref{G-pr} for every $s\in I_L$
\beq\label{M}
M(s):=G_-(s)^{-1}(G(s)^\top)^{-1}=2\left(A(s)^\top+\bI_N\right)^{-1}\,,
\eeq
the matrix $A(s)$ being given by \eqref{As}, so that $M(s)\neq \bI_N$ if and only if $s\in S(G)$.
For that reason, the problem in \eqref{sist-wrong} is well-defined if we impose that the matrix $A(s)$ satisfies property \eqref{detAs}
at every Jump point $s\in S(G)$.

Notice that the latter assumption is consistent with the classical result by Cayley reported in Proposition \ref{P-Cayley}.
However, we now see that even in the planar case the system \eqref{sist-wrong} fails to give the right approach to deal with
a datum $\O(s)$ with jump points.

For simplicity of notation, from now one we shall assume $\bG=\bI_N$ in \eqref{eqND}.
The more general case is trivially obtained by means of a right multiplication of the solution we obtain with the initial datum $\bG\in SO(N)$.
%
%
\subsection{The 2D case}
In low dimension $N=2$, any function $\O \in BV( I_L, Sk(2))$ is given by
\begin{equation}
\label{Omega2}
\O(s)=\theta(s)\,\wih\bJ\,,\quad \wih\bJ:=\left(
     \begin{array}{cc}
       0 & 1 \\
       -1 & 0 \\
     \end{array}
   \right)
\end{equation}
for some $\theta \in BV(I_L)$. In order that solutions $G\in BV(I_L,SO(2))$ satisfy eq. \eqref{detAs} at every point $s\in S(G)$, so that problem \eqref{eqND} makes sense, we \emph{assume that $\t$ is strictly increasing in $I_L$ and that for every discontinuity point $s \in S(\t)$ of $\theta$}
\beq
\label{small-jumps}
\theta(s^+)-\theta(s^-)<\pi\,,\qquad \theta(s^\pm):=\lim_{\s\to s^\pm}\theta(\s)\,.
\end{equation}

Therefore, letting $\gt(s)=G(s){\bf e}_1$ and $\gn(s)=G(s){\bf e}_2$, where $({\bf e}_1,{\bf e}_2)$ denote the (column) vectors of the canonical orthonormal basis in $\gR^2$, problem \eqref{eqND},
where $\bG=\bI_2$, is equivalent to
\[
\left\{ \begin{array}{l}
D\gt =\gn\,D\theta \\
D\gn=-\gt\,D\theta \\
\gt(0)={\bf e}_1 \,,\quad \gn(0)={\bf e}_2 \,.
\end{array} \right.
\]
If $\theta \in W^{1,1}(I_L)$, letting  $k(s):=\theta'(s)$ we obtain the system
\[
\left\{ \begin{array}{l}
\gt'(s) =k(s)\,\gn(s) \\
\gn'(s)=-k(s)\,\gt(s) \\
\gt(0)={\bf e}_1 \,,\quad \gn(0)={\bf e}_2
\end{array} \right.
\]
that has to be intended in the $\calL^1$-almost everywhere sense.
%
%
%
In general, we have:
\bp\label{P2D} Assume in addition that $\t$ is continuous on $I_L$.
Then the problem \eqref{eqND}, where $\bG=\bI_2$,
has a unique solution in the class of continuous maps $G\in BV(I_L,\gR^{2\tim 2})$, given by
\beq\label{sol2D}
G(s):=\left(\begin{array}{cc}
       \cos \a(s) & - \sin\a(s) \\
       \sin\a(s) & \cos\a(s) \\
     \end{array}
   \right) \in BV(I_L,SO(2))\,,
\eeq
where $\a(s)=\theta(s)-\theta(0)$.
Moreover, letting $\gt(s)= G(s) {\bf e_1}$, there exists a unique (up to rigid motions) rectifiable curve $\gc:\bar I_L\to\gR^2$ parameterized in arc-length, with tantrix $\gc'=\gt$ and with finite total curvature
\[
\TC(\gc)= |D\gt|(I_L)=\t(L)-\t(0)\,.
\]
\ep
\bpf We have seen that any solution $G\in BV(I_L,\gR^{2\tim 2})$ 
to the system \eqref{eqND} is a continuous map in $BV(I_L,SO(2))$,
and hence it can be written as in \eqref{sol2D} for some continuous angle function $\a\in BV(I_L)$.
Therefore, since problem \eqref{eqND} is equivalent to \eqref{eqND-cont}, a direct computation yields that $D\a=D\t$,
whence $\a(s)=\theta(s)+c$, and the initial condition gives $c=-\t(0)$.
The other assertions readily follow.
\epf

Now, adding to the angle function $\t(s)$ an isolated jump of small size, i.e. by letting
$$ \wid\t(s)=\theta(s)+d\cdot\chi_{(s_0,L)}(s) $$
for some $s_0\in I_L$ and $0<d<\pi$, where $\chi_I$ denotes the characteristic function of a set $I\sb I_L$, we have $D\wid\t=D\theta+d\,\d_{s_0}$. Therefore the solution $\wid G(s)$ satisfies $\wid G(s)=G(s)$ for $s\in(0,s_0)$, where $G(s)$ is given by \eqref{sol2D} with $\a(s)=\theta(s)-\theta(0)$.

If we consider the system \eqref{sist-wrong}, using the notation from Example \ref{E2D}, and on account of \eqref{M}, we obtain the equation
$$ \left(\begin{array}{cc} 0 & \sin d \\ -\sin d & 0 \\\end{array}\right) \d_{s_0}
=
\left(\begin{array}{cc}
1 & -\sin d/(1+\cos d) \\
\sin d/(1+\cos d) & 1 \\\end{array}\right)
\left(\begin{array}{cc} 0 & -d \\ d & 0 \\ \end{array}\right)
\d_{s_0}
$$
for the Jump part of the derivative, that does not make sense.

In fact, in the planar case $N=2$, a Jump of size $d$ for the angle function at $s_0$ corresponds to a rotation matrix $R(d)\in SO(2)$ of angle $d$. Therefore, the solution $\wid G(s)$ has a Jump point at $s_0$, and its prolongation to $(s_0,L)$ is given by the unique solution to the system
$$
\left\{ \begin{array}{l}
DG= - G\,D\O \\
G(s_0^+)=R(d)\,G(s_0^-)\,,
\end{array} \right. \quad s\in (s_0,L)\,,
$$
i.e. by the function $\wid G(s)$ in \eqref{sol2D}, where this time $\a(s)=\theta(s)-\theta(0)+ d\cdot\chi_{(s_0,L)}(s)$.
Finally, notice that
$$ D\wid G = \nabla \wid G\,\calL^1+ (R(d)-\bI_2)\,\wid G(s_0^-)\,\d_{s_0}\,, $$
where $\nabla \wid G=\nabla G$ for $\calL^1$-a.e. $s\in(0,s_0)$, and $\nabla \wid G=R(d)\nabla G$ for $\calL^1$-a.e. $s\in(s_0,L)$.

\smallskip

Now, denoting for simplicity $G_-=G(s_0^-)$ and $R=R(d)$, we have $G(s_0^+)=RG_-$, the precise representative is
$\wid G(s_0)=2^{-1}(R+\bI_2)G_-$ and the Jump component of the derivative is $D^J\wid G=(R-\bI_2)G_-\d_{s_0}$, so that by the commutativity of the product of matrices in $SO(2)$ we readily obtain that
$$ -\wid G^\top D^J\wid G=\sin d\,\wih\bJ\,\d_{s_0}\,, $$
where $\wih\bJ$ is defined in \eqref{Omega2}.
Therefore, the atomic component of the equation between measures in \eqref{eqND} makes sense provided that we replace $D^J\O$ with the atomic measure
$ D^{at}\O:=\sin d\,\wih\bJ\,\d_{s_0}$.
%
%
Notice in fact that according to Example \ref{E2D}, if $R(\DD\a)\in SO(2)$ corresponds to a planar rotation with angle $\DD\a$, then
$$ \frac 12\left(R(\DD\a)-R(\DD\a)^\top\right) = \sin \DD\a\,\wih\bJ =
\left(\begin{array}{cc} 0 & \sin \DD\a \\
-\sin \DD\a & 0 \\  \end{array} \right)\,.$$

In a similar way, if $D^J\t=\sum_{i=1}^n d_i\,\d_{s_i}$, where the discontinuity points $\{s_i\}_{i=1}^n\sb I_L$ are pairwise disjoint, and $0<d_i<\pi$ for $i=1,\ldots,n$, at the right-hand side of eq. \eqref{eqND}
we have to replace the Jump component $D^J\O$ of the datum in \eqref{Omega2} with the atomic measure
\beq\label{atom-2D}
D^{at}\O:=\sum_{i=1}^n \sin d_i\,\wih\bJ\,\d_{s_i}\,.
\eeq
\par
The general case, i.e. when the Jump of the datum $\t$ is an infinite (countable) set, can be treated in a way similar to the one in the physical dimension $N=3$ that we shall discuss below.
%
%
\subsection{A case study with a Jump point in 3D}
%
%
We now analyze in detail an example in dimension $N=3$ where both the functions $\t$ and $\f$ in \eqref{Omega3} are smooth except at a singular point, where at least one of them has a Jump.
Our purpose is to understand a way to make sense of eq. \eqref{eqND}, when the Jump component $D^J\O$ of the datum is non-zero.

We thus let $d\geq 0$ and $\ttt\in\gR$ be such that 
\beq\label{dtsmall}
0<\sqrt{d^2+\ttt^2}<\pi\,. \eeq
As we shall see below, the bound in \eqref{dtsmall} allows us to recover the invertibility property of the solution $G(s)$ at the Jump point $s=0$.
We define
$$ \t(s)=\left\{\ba{ll} s & \textrm{if } s<0  \\ s+d & \textrm{if } s>0\,,\ea\right.
\qquad
\f(s)=\left\{\ba{ll} 0 & \textrm{if } s<0  \\ \ttt & \textrm{if } s>0\,,\ea\right.
$$
so that we have $D\t\pri I=\calL^1\pri I+d\,\d_0$ and $D\f\pri I=\ttt\,\d_0$, where this time we denote $I=(-1,1)$.

We shall construct a rectifiable curve $\gc:I\to\gR^3$ parameterized in
arc-length by the limit of solutions of the Frenet--Serret system
associated to the derivatives of smooth approximations of the functions $\t$ and $\f$.
Moreover, the total curvature and the total absolute torsion of $\gc$ are linked to the distributional derivatives of the vector valued function $(\t,\f):I\to\gR^2$.

More precisely, setting
$\gt(s):=\gc'(s)$ and $\gn(s)=\gt'(s)/|\gt'(s)|$, for $s\neq 0$, we have
$$ |D\gn|(I)=|D(\t,\f)|(I)=2+\sqrt{d^2+\ttt^2}\,.
$$
However, it turns out that the singular part of the distributional derivative of the tantrix $\gt$, as well as the one of the binormal vector $\gb:=\gt\tim\gn$, depends on the Jump of the couple $(\t,\f)$.

For $\e>0$ small, we define $\t_\e,\f_\e:\ol I\to\gR$ by the Lipschitz-continuous average functions
\beq\label{approx}
\t_\e(s):=\frac 1{2\e}\int_{s-\e}^{s+\e}\t(\s)\,\dd\s\,,\quad \f_\e(s):=\frac 1{2\e}\int_{s-\e}^{s+\e}\f(\s)\,\dd\s\,,
\eeq
so that we have
$$ \t_\e(s)=\left\{\ba{ll} s & \textrm{if } s\leq -\e  \\
\ds\frac{2\e+d}{2\e}\,s+\frac d 2 & \textrm{if } |s|\leq \e \\
s+d & \textrm{if } s\geq \e\,,\ea\right.
\qquad
\f_\e(s)=\left\{\ba{ll} 0 & \textrm{if } s\leq -\e  \\
\ds\frac \ttt{2\e}\,s+\frac \ttt 2 & \textrm{if } |s|\leq \e \\
\ttt & \textrm{if } s\geq \e\,. \ea\right.
$$

We now consider the Cauchy problem
\beq\label{CP-example}
\left\{ \begin{array}{l}
G'(s) = - (G^\top(s))^{-1} \O'_\e(s) \\
G(0)=\bI_3\,,
\end{array}\right.
\quad s\in \ol I\,,
\eeq
where $\O_\e:I\to Sk(3)$ is given by \eqref{Omega3} with $\t=\t_\e$ and $\f=\f_\e$, so that
for $s\neq \pm\e$ we have
$$
\O'_\e(s):=\left(\begin{array}{ccc}
0 & \t'_\e(s) & 0 \\ -\t'_\e(s) & 0 & \f'_\e(s) \\ 0  & -\f'_\e(s) & 0  \\
 \end{array} \right) $$
where
$$
\t'_\e(s)=\left\{\ba{ll}
\ds\frac{2\e+d}{2\e} & \textrm{if } |s| < \e \\
1 & \textrm{if } \e<|s|\leq 1\,,\ea\right.
\quad
\f'_\e(s)=\left\{\ba{ll}
\ds\frac \ttt{2\e} & \textrm{if } |s| < \e \\
0 & \textrm{if } \e<|s|\leq 1\,. \ea\right.
$$

For $s\in(-\e,\e)$, according to the computation in Appendix \ref{App:B}, using the Rodrigues formula
we find the unique solution $G_\e$ given by the exponential map
$$
G_\e(s)=\left(\begin{array}{ccc}
c_\e(s)+{\tau_\e}^2(1- c_\e(s)) & k_\e\, s_\e(s) & k_\e\,\tau_\e\,(1-c_\e(s) ) \\
-k_\e\,s_\e(s) & c_\e(s) & \tau_\e\,s_\e(s) \\
k_\e\,\tau_\e\,(1-c_\e(s))  & -\tau_\e\,s_\e(s) & c_\e(s)+{k_\e}^2(1-c_\e(s))  \\
 \end{array} \right) \,, $$
where for the sake of brevity we have denoted $c_\e(s):=\cos ((\l_\e/\e) s)$,  $s_\e(s):=\sin ((\l_\e/\e) s)$, and
$$ \l_\e:=\frac{\sqrt{(2\e+d)^2+\ttt^2}}{2}\,,\quad k_\e:=\ds\frac {2\e+d}{2\l_\e}\,,\quad
\tau_\e:= \frac {\ttt}{2\l_\e}\,.
$$
In particular, at times $s=\pm\e$ we have
$$
G_\e(\pm\e)=\left(\begin{array}{ccc}
\cos \l_\e+{\tau_\e}^2(1-\cos \l_\e) & \pm k_\e\sin \l_\e  & k_\e\,\tau_\e\,(1-\cos \l_\e) \\  \mp k_\e\,\sin \l_\e & \cos \l_\e & \pm\tau_\e\,\sin \l_\e \\  k_\e\,\tau_\e\,(1-\cos \l_\e)  & \mp\tau_\e\,\sin \l_\e & \cos \l_\e+{k_\e}^2(1-\cos \l_\e)  \\
 \end{array} \right)\,.
$$
Therefore, the solution $G_\e(s)$ extends to the interval $\ol I$ by letting
%
$$ G_\e(s)=\ds \left(\begin{array}{ccc}
\cos (s\mp\e) & -\sin (s\mp\e) & 0 \\  \sin (s\mp\e) & \cos (s\mp\e) & 0 \\  0 & 0 & 1  \\
 \end{array} \right) G_\e(\pm\e) \quad \textrm{if } \pm s > \e\,.
$$
\par
By taking the limit as $\e\to 0^+$, we thus obtain the function
$$ G(s):=\ds \left(\begin{array}{ccc}
\cos s & -\sin s & 0 \\  \sin s & \cos s & 0 \\  0 & 0 & 1  \\
 \end{array} \right) G(0^\pm) \quad \textrm{if } 0< \pm s \leq 1\,,
$$
where we have denoted
$$ G(0^\pm):=\lim_{\e\to 0^+} G_\e(\pm\e)\,. $$
\par Finally, setting
$$ \gc_\e(s):=0_{\gR^3}+\int_0^s\gt_\e(\s)\,\dd \s\,,\quad \gt_\e(s):=G_\e(s)\,{\bf e}_1\,,
$$
it turns out that $\gc_\e$ strongly converges to the expected solution curve $\gc$.
\subsection{Analysis at the Jump point}
The function $G(s)$ this way obtained is continuous at $s\neq 0$, whereas its right and left limits at $s=0$ are
$$
G(0^\pm)=\left(\begin{array}{ccc}
\cos\a+\frac {\ttt^2}{d^2+\ttt^2}\,(1-\cos \a) & \pm \frac d{\sqrt {d^2+\ttt^2}}\sin \a  & \frac {d\ttt}{d^2+\ttt^2}\, (1-\cos \a) \\
 \mp \frac d{\sqrt {d^2+\ttt^2}}\,\sin \a & \cos \a & \pm \frac \ttt{\sqrt {d^2+\ttt^2}}\,\sin \a \\
 \frac {d\ttt}{d^2+\ttt^2}\, (1-\cos \a)  & \mp\frac \ttt{\sqrt {d^2+1}}\,\sin \a & \cos \a+\frac {d^2}{d^2+\ttt^2}\,(1-\cos\a)  \\
 \end{array} \right)\,, \quad \a:=\frac{\sqrt{d^2+\ttt^2}}2\,, $$ 
that correspond to a rotation of angle $\mp\a$ around the axis oriented by $(\ttt/\sqrt{d^2+\ttt^2},0,d/\sqrt{d^2+\ttt^2})$.

According to \eqref{As}, the matrix $A(0):=G(0^-)^\top G(0^+)$ is a rotation of angle $-2\a$. 
Since the bound in \eqref{dtsmall} gives $ |2\a|<\pi$,
we infer that eq. \eqref{detAs} is satisfied.
Therefore, the precise representative of $G(s)$ satisfies
$$ G(0)=\frac 12\,
G(0^-) \left(A(0)+\bI_3\right)
$$
and hence $G(0)$ is an invertible matrix. Moreover, the Jump part of the derivative is given by
$$ D^J G =
G(0^-) \left(A(0)-\bI_3\right)
\,\d_0\,. $$
%
We thus compute
$$ G^\top D^JG = \frac 12\,\left(A(0)^\top+\bI_3\right)G(0^-)^\top G(0^-)\left(A(0)-\bI_3\right)
\,\d_0
= -\frac 12 \left(A(0)^\top-A(0)\right)\,\d_0 \,.$$
Therefore, the Jump part of equation \eqref{eqND} makes sense provided that we replace the Jump component $D^J\O$ with the atomic measure
$$ D^{at}\O:= \frac 12 \left(A(0)^\top-A(0)\right)\,\d_0\,. $$
\par Notice that the datum $\O(s)$ in \eqref{Omega3} can be written in terms of the generators in \eqref{matrix-J} by the formula
$$ \O(s)=-\t(s)\,\bJ_3-\f(s)\,\bJ_1\,. $$
Since moreover the matrix $A(0)^\top$ is a rotation of angle $2\a=\sqrt{d^2+\ttt^2}$ around the axis oriented by the unit vector
$(\ttt/\sqrt{d^2+\ttt^2},0,d/\sqrt{d^2+\ttt^2})$, we correspondingly obtain
\beq\label{atom-3D}
D^{at}\O=- \left( \frac {d\,\sin\sqrt{d^2+\ttt^2}}{\sqrt {d^2+\ttt^2}} \,\bJ_3 + \frac {\ttt\,\sin\sqrt{d^2+\ttt^2}}{\sqrt {d^2+\ttt^2}}\,\bJ_1\right)\,\d_0\,.
\eeq
In particular, if $\ttt=0$ or $d=0$, we respectively obtain
$$ D^{at}\O=-\sin d\,\bJ_3\,\d_0\,,\quad D^{at}\O=-\sin \ttt\,\bJ_1\,\d_0\,,$$
where the first formula may be compared with eq. \eqref{atom-2D} for the planar case.
\par
When $\ttt=0$ and $0<d<\pi$, we have
$$
G(0^\pm)=\left( \gt^\pm \vert \gn^\pm \vert \gb^\pm \right):=\left(\begin{array}{ccc}
\cos(d/2) & \pm\sin(d/2) & 0 \\
 \mp\sin(d/2) & \cos(d/2) & 0 \\
 0 & 0  & 1  \\
 \end{array} \right)
$$
and the angle between $\gt^\pm$ is equal to $d$.
This case corresponds to a planar curve $\gc$ with a corner point at $s=0$ of turning angle $d$, and curvature equal to one elsewhere.

When $d=0$ and $0<|\ttt|<\pi$, instead, we have
$$
G(0^\pm)=\left( \gt^\pm \vert \gn^\pm \vert \gb^\pm \right):=\left(\begin{array}{ccc}
1 & 0 & 0 \\
 0 & \cos(\ttt/2) & \pm\sin(\ttt/2) \\
 0 & \mp\sin(\ttt/2)  & \cos(\ttt/2)  \\
 \end{array} \right)\,.
$$
The solution curve $\gc$ is again piecewise smooth, with zero torsion and curvature equal to one for $s\neq 0$.
Moreover, it is given by two circular arcs of radius one meeting at the point $\gc(0)$ but belonging to two different planes,
the first one normal to the vector $\gb^-:=(0,-\sin(\ttt/2),\cos(\ttt/2))$ and the second one to the vector $\gb^+:=(0,\sin(\ttt/2),\cos(\ttt/2))$.
Since $\gt^-=\gt^+$, no curvature appears at the point $\gc(0)$, whereas the jumps $\gn^+-\gn^-$ and $\gb^+-\gb^-$ are due to the occurrence of a jump in the function $\f$ at $s=0$.
More precisely, both the angles between $\gb^\pm$ and $\gn^\pm$ are equal to $|\ttt|$, whereas the corresponding Euclidean distances are equal to
$2\sin(|\ttt|/2)$.

More generally, when both $d$ and $\ttt$ are non-zero, the angle between $\gt^\pm$ is
$$ \widehat{\gt^+\vert\gt^-}=\arccos\left(\frac{\ttt^2+d^2\cos\sqrt{d^2+\ttt^2}}{d^2+\ttt^2}\right)\,, $$
the angle between $\gn^\pm$ is
$ \widehat{\gn^+\vert\gn^-}=\sqrt{d^2+\ttt^2} $,
and the angle between $\gb^\pm$ is
$$ \widehat{\gb^+\vert\gb^-}=\arccos\left(\frac{d^2+\ttt^2\cos\sqrt{d^2+\ttt^2}}{d^2+\ttt^2}\right)\,. $$
In that case, either the turning angle at the corner point $\gc(0)$ or the distance between the planes containing the two circular arcs of $\gc$ meeting at $\gc(0)$, both depend on $d$ and $\ttt$.
These features are consistent with our choice in formula \eqref{atom-3D}.

We now see that both the angles between $\gt^\pm$ and 
$\gb^\pm$ are bounded by the angle between $\gn^\pm$. 
\bp\label{Pangle} With the previous assumptions, we have:
$$ 0\leq\widehat{\gt^+\vert\gt^-}\leq \sqrt{d^2+\ttt^2}\quad\text{and}\quad  
0\leq\widehat{\gb^+\vert\gb^-}\leq \sqrt{d^2+\ttt^2}\,.
$$
\ep
\bpf Using that
$$ 1=\frac{\ttt^2+d^2}{d^2+\ttt^2}\geq \frac{\ttt^2+d^2\cos\sqrt{d^2+\ttt^2}}{d^2+\ttt^2}\geq \frac{(\ttt^2+d^2)\,\cos\sqrt{d^2+\ttt^2}}{d^2+\ttt^2}= \cos\sqrt{d^2+\ttt^2}$$
we obtain
$$ 0\leq \arccos\left(\frac{\ttt^2+d^2\cos\sqrt{d^2+\ttt^2}}{d^2+\ttt^2}\right)\leq \sqrt{d^2+\ttt^2}\,, $$
that gives the first inequality. The second one is proved in a similar way.
\epf

Finally, by the previous observations we readily obtain the explicit formulas for the total curvature and absolute torsion of the solution curve, namely:
$$ \ba{l} \ds \TC(\gc)=2+ \arccos\left(\frac{\ttt^2+d^2\cos\sqrt{d^2+\ttt^2}}{d^2+\ttt^2}\right)\,, \\
\ds \TAT(\gc)=\arccos\left(\frac{d^2+\ttt^2\cos\sqrt{d^2+\ttt^2}}{d^2+\ttt^2}\right)\,. 
\ea $$
Notice that if $\pi/2\leq\sqrt{d^2+\ttt^2}<\pi$ and $d^2+\ttt^2\cos\sqrt{d^2+\ttt^2}<0$, the contribution to the total absolute torsion is equal to $\pi-\widehat{\gb^+\vert\gb^-}$, that is the distance between $\gb^\pm$ when embedded in the projective plane or, equivalently, the distance between the (unoriented) osculating planes generated by $(\gt^+,\gn^+)$
and $(\gt^-,\gn^-)$.
%
%
\section{A more general existence theorem}\label{Sec:gen}
In this section we deal with the general case when the datum $\O\in BV(I_L,Sk(3))$ fails to be continuous. 
We thus let $N=3$, and the functions $\t,\f$ in \eqref{Omega3} are assumed to be
in $BV(I_L)$, with $\theta(s)$ strictly increasing.
We first deal with the case when the datum $\O(s)$ has a finite number of discontinuity points, and we then discuss the general case. Finally, for simplicity we assume $\bG=\bI_3$ in \eqref{eqND}.
\subsection{The case of a finite discontinuity set}
Assume that the Jump set $S(\O)$ is finite, namely
$$ S(\O)=\{0<s_1<s_2<\ldots<s_m<L\}\,. $$
We recall the notation $\bar D u$ for the diffuse component of the measure derivative, so that
$$ \bar D\O:=D^a\O+D^C\O\,,\quad \bar DG:=D^aG+D^CG\,. $$
The Jump component of $D\O$ is given by
$$ D^J\O=\sum_{k=1}^m [\O(s_k)]\,\d_{s_k}\,,\quad [\O(s_k)]=\O(s_k^+)-\O(s_k^-)\,,\quad k=1,\ldots,m\,. $$
Therefore, the $BV$ functions $\t(s)$ and $\f(s)$ can be discontinuous only at the points in $S(\O)$, in such a way that $S(\O)=S(\t)\cup S(\vf)$. Following the example from the previous section, for every $k$ we can write
$$ [\O(s_k)]=
\left(\begin{array}{ccc}
0 & d_k & 0 \\ -d_k & 0 & \ttt_k \\ 0  & -\ttt_k & 0  \\
 \end{array} \right)\,,\quad d_k:=[\theta(s_k)]\,,\quad \ttt_k:=[\phi(s_k)]\,.
$$
Therefore, $d_k=0$ or $\ttt_k=0$ if $\t(s)$ or $\f(s)$ is continuous at $s=s_k$.
Finally, we assume that the Jumps of the vector function $(\t,\f)$ are smaller that $\pi$, so that definitely
\beq\label{smalljumps}
0<\sqrt{ {d_k}^2+{\ttt_k}^2}<\pi\quad\fa\,k=1,\ldots,m\,.
\eeq

According to \eqref{atom-3D}, we introduce the atomic measure
$$
D^{at}\O:=- \sum_{k=1}^m\bB_k \,\d_{s_k}
$$
where, using the notation in \eqref{matrix-J} for the generators $\bJ_1$, $\bJ_3$, we set
$$\bB_k:=\left( \frac {d_k\,\sin\sqrt{{d_k}^2+{\ttt_k}^2}}{\sqrt {{d_k}^2+{\ttt_k}^2}} \,\bJ_3 + \frac {\ttt_k\,\sin\sqrt{{d_k}^2+{\ttt_k}^2}}{\sqrt {{d_k}^2+{\ttt_k}^2}}\,\bJ_1\right)\,,\quad k=1,\ldots,m\,.
$$
Recall that if $\ttt_k=0$ or $d_k=0$ for some $k$, we respectively have
$$ \bB_k=\sin d_k\,\bJ_3 \quad \textrm{or}\quad \bB_k=\sin \ttt_k\,\bJ_1\,.  $$
Finally, we introduce the $Sk(3)$-valued Borel measure $\wid D\O$ given by
$$ \wid D\O:= \bar D\O+D^{at}\O\,, $$
and in the sequel we identify a function $G\in BV(I_L,SO(3))$ with its precise representative, compare \eqref{G-pr}.
\bt\label{T-disc} Under the previous notation, if \eqref{smalljumps} holds, then
the system
$$
\left\{ \begin{array}{l}
DG= - (G^\top)^{-1} \wid D\O \\
G(0)=\bI_3\,,
\end{array} \right.
$$
has a unique solution $G\in BV(I_L,\gR^{3\tim 3})$ such that $G(s)\in SO(3)$ for every $s\in I_L\sm S(\O)$.
\et
\bpf 
For $k=1,\ldots,m$, denote
$$ 2\a_k:=\sqrt{{d_k}^2+{\ttt_k}^2} \,, \quad \gv_k:=(\ttt_k/(2\a_k),0,{d_k}/(2\a_k))^\top\,. $$
The function $G(s)$ on the first interval $(0,s_1)$ of continuity of $\O$ is given by applying Theorem \ref{Tmain}.
Let $\bG_1\in SO(3)$ denote the rotation matrix of angle $2\a_1$ around the axis oriented by
$G(s_1^-)\gv_1$. The function $G(s)$ on the second interval $(s_1,s_2)$ is then determined by solving the problem
$$
\left\{ \begin{array}{l}
DG= - (G^\top)^{-1} D\O \\
G(s_1^+)=\bG_1\,G(s_1^-)\,,
\end{array} \right. \quad s\in(s_1,s_2)\,.
$$
Let now $\bG_2\in SO(3)$ denote the rotation matrix of angle $2\a_2$ around the axis oriented by
$G(s_2^-)\gv_2$.
The function $G(s)$ on the third interval $(s_1,s_2)$ is determined by solving the problem
$$
\left\{ \begin{array}{l}
DG= - (G^\top)^{-1} D\O \\
G(s_2^+)=\bG_2\,G(s_2^-)\,,
\end{array} \right. \quad s\in(s_2,s_3)\,.
$$
Proceeding this way, we readily find a solution to the given problem. In fact, for $k=1,\ldots,m$, letting
$A(s_k):=G(s_k^-)^\top G(s_k^+)$, we have
$$ G(s_k)=\frac 12\,
G(s_k^-) \left(A(s_k)+\bI_3\right)\,.
$$
As in the previous case study, property \eqref{smalljumps} implies that eq. \eqref{detAs} is satisfied by $A(s_k)$ for every $k$, so that the matrices $G(s_k)$ are invertible. 
Moreover, the Jump part of the derivative is given by
$$ D^J G =\sum_{k=1}^m
G(s_k^-) \left(A(s_k)-\bI_3\right)
\,\d_{s_k}\,. $$
We thus have
$$ G^\top D^JG = -\frac 12\,\sum_{k=1}^m
\left(A(s_k)^\top-A(s_k)\right)\,\d_{s_k}\,,$$
where for $k=1,\ldots,m$
$$\frac 12\,
\left(A(s_k)^\top-A(s_k)\right)=\bB_k\,. $$
Therefore, we obtain the identity between the diffuse and atomic components, namely:
$$ \bar DG = - (G^\top)^{-1} \bar D\O\quad\textrm{and}\quad D^JG= - (G^\top)^{-1} D^{at}\O\,. $$
Uniqueness readily follows on account of Theorem \ref{Tmain}.
\epf
\par Recalling the notation in \eqref{tnb}, we now define
$$ \gc(s):=0_{\gR^3}+\int_0^s\gt(\s)\,d\s\,,\quad s\in \ol I_L\,. $$
Then $\gc:\ol I_L\to\gR^3$ is a piecewise $C^1$ rectifiable curve parameterized in arc-length and with tantrix $\gc'(s)=\gt(s)$ for every $s\in I_L\sm S(\O)$, i.e. when $s\neq s_k$ for every $k$.

By our construction, for $k=1,\ldots,m$, the angle between $\gt^\pm(s_k)$ is
$$ \widehat{\gt^+(s_k)\vert\gt^-(s_k)}=\arccos\left(\frac{{\ttt_k}^2+{d_k}^2\cos\sqrt{{d_k}^2+{\ttt_k}^2}}{{d_k}^2+{\ttt_k}^2}\right)\,, $$
the angle between $\gn^\pm(s_k)$ is
$ \widehat{\gn^+(s_k)\vert\gn^-(s_k)}=\sqrt{{d_k}^2+{\ttt_k}^2} $,
and the angle between $\gb^\pm(s_k)$ is
$$ \widehat{\gb^+(s_k)\vert\gb^-(s_k)}=\arccos\left(\frac{{d_k}^2+{\ttt_k}^2\cos\sqrt{{d_k}^2+{\ttt_k}^2}}{{d_k}^2+{\ttt_k}^2}\right)\,. $$
Finally, recall the mutual decomposition into diffuse and atomic components of the derivatives of the column vectors of $G$, so that
$$ D\gt=\bar D\gt+D^J\gt\,,\quad D\gb=\bar D\gb+D^J\gb \,. $$

By applying separately Theorems \ref{Tmain} and \ref{T-tor} on each component of $I_L\sm S(\O)$, and arguing as in the example from the previous section, we thus readily obtain the following
\bp\label{P-disc}
The solution curve $\gc$ has finite total curvature given by
$$ \TC(\gc)=|\bar D\gt|(I_L)+\sum_{k=1}^m \arccos\left(\frac{{\ttt_k}^2+{d_k}^2\cos\sqrt{{d_k}^2+{\ttt_k}^2}}{{d_k}^2+{\ttt_k}^2}\right)
$$
and finite total absolute torsion equal to
$$ \TAT(\gc)=|\bar D\gb|(I_L)+\sum_{k=1}^m \arccos\left(\frac{{d_k}^2+{\ttt_k}^2\cos\sqrt{{d_k}^2+{\ttt_k}^2}}{{d_k}^2+{\ttt_k}^2}\right)\,.
$$
\ep
\br\label{Rapprox} For $\e>0$ small, we can define $\t_\e,\f_\e:\ol I\to\gR$ by the Lipschitz-continuous average functions as in \eqref{approx} and solve the corresponding Cauchy problem \eqref{CP-example}. Since $\O$ is continuous outside a finite set, arguing as in the case study from the previous section it is not difficult to show that the solution curve $\gc_\e$ converges (along a sequence $\e_h\searrow 0$) to the curve $\gc$ previously obtained. \er
\par For future use, we finally point out the following estimates:
\bc\label{Cangle} We have:
$$ \TC(\gc)\leq|\bar D\gt|(I_L)+2^{-1/2}|D^J\O|(I_L)
\quad\text{and}\quad\TAT(\gc)\leq|\bar D\gb|(I_L)+2^{-1/2}|D^J\O|(I_L)\,.
$$
\ec
\bpf Arguing as in the proof of Proposition \ref{Pangle} we obtain the bound
$$ 0\leq \arccos\left(\frac{{\ttt_k}^2+{d_k}^2\cos\sqrt{{d_k}^2+{\ttt_k}^2}}{{d_k}^2+{\ttt_k}^2}\right)\leq \sqrt{{d_k}^2+{\ttt_k}^2} $$
for $k=1,\ldots,m$. Thefefore, the first inequality follows since by the previous notation
we have
$$|D^J\O|=\sqrt 2\,\sum_{k=1}^m\sqrt{{d_k}^2+{\ttt_k}^2}\,.
$$
The second inequality is proved in a similar way. 
\epf
\subsection{The general case}
We finally remove the assumption that $\O(s)$ has a finite Jump set. As in \eqref{smalljumps}, for every $s\in S(\O)$ we require
$$ 0<\sqrt{d(s)^2+\ttt(s)^2}<\pi\,,\quad d(s)=[\t(s)],\quad \ttt(s)=[\f(s)]\,.
$$
According to the notation in \eqref{matrix-J}, we define the atomic measure
$$
D^{at}\O:=- \sum_{s\in S(\O)} \left( \frac {d(s)\,\sin\sqrt{d(s)^2+\ttt(s)^2}}{\sqrt {d(s)^2+\ttt(s)^2}} \,\bJ_3 + \frac {\ttt(s)\,\sin\sqrt{d(s)^2+\ttt(s)^2}}{\sqrt {d(s)^2+\ttt(s)^2}}\,\bJ_1\right) \,\d_{s}\,.
$$
We then proceed by approximation, and for every $n\in\Nat^+$ we let
$$ S_n(\O):=\{s\in S(\O) \,:\, |D^J\O|(s)>1/n\}\,. $$
Property $|D^J\O|<\i$ implies that $S_n(\O)$ is a finite set for every $n$.
Moreover, $\{S_n(\O)\}$ is an increasing sequence converging to $S(\O)$,
so that $|D^J\O|(I_L\sm S_n(\O))\to 0$ as $n\to\i$.
We then define
$$ D^{at}_n\O:=(D^{at}\O)\pri S_n(\O)\,,\quad n\in\Nat^+ $$
and consider the approximate problem
\beq\label{appr-pb}
\left\{ \begin{array}{l}
DG= - (G^\top)^{-1} \left(\bar D\O+D^{at}_n\O\right) \\
G(0)=\bI_3\,,
\end{array} \right. \quad n\in\Nat^+\,.
\eeq
\par By Theorem \ref{T-disc}, we infer that the latter system has a unique solution
$G_n\in BV(I_L,\gR^{3\tim 3})$ such that $G_n(s)\in SO(3)$ for every $s\in I_L\sm S_n(\O)$. Denoting by $\gt_n$, $\gn_n$, and $\gb_n$ the ordered column vectors of the matrix $G_n(s)$, for every $n$ we define $\gc_n:\ol I_L\to\gR^3$ by
\beq\label{cn}
\gc_n(s):=0_{\gR^3}+\int_0^s\gt_n(\s)\,d\s\,,\quad s\in \ol I_L\,,
\eeq
so that Proposition \ref{P-disc} applies to every piecewise $C^1$ curve $\gc_n$.
\bp\label{P-cn} The sequence $\{\gc_n\}$ converges uniformly to a rectifiable curve $\gc$ of length $L$ and with finite total curvature and total absolute torsion.
\ep
\bpf We have $||DG_n|(I_L)-|DG_{n+1}|(I_L)|\to 0$ as $n\to \i$. Therefore, $\{\gc_n\}$ is a Cauchy sequence w.r.t. the Frech\'et distance, and hence it converges to a rectifiable curve $\gc$ of length $L$. 
Furthermore, by Corollary \ref{Cangle} we obtain the estimate
$$ \TC(\gc_n)\leq |\bar D\gt|(I_L)+2^{-1/2}|D^J\O|(I_L)
\quad\text{and}\quad\TAT(\gc_n)\leq|\bar D\gb|(I_L)+2^{-1/2}|D^J\O|(I_L)
$$
for every $n\in\Nat^+$, whence
$$ \sup_n\left(\TC(\gc_n)+\TAT(\gc_n)\right)<\i\,. $$
Therefore, the second assertion follows by the lower semicontinuity of the total curvature and total absolute torsion w.r.t. the convergence in the Frech\'et distance.
\epf

In principle, the curve $\gc$ depends on the choice of the approximating sequence $\{\gc_n\}$ in \eqref{cn}. Assume now that $\{\hat\gc_h\}$ is another sequence of curves obtained as above, but this time w.r.t. the solutions $\{\hat G_h(s)\}$ to some approximate problems as in \eqref{appr-pb}, where we replace the term  
$ D^{at}_n\O$ by
$$D^{at}_h\O=(D^{at}\O)\pri S_h(\O)\,,\quad h\in\Nat^+ $$
and $\{S_h(\O)\}$ is any sequence satisfying the following properties:
\begin{itemize}
\item 
$S_h(\O)$ is a finite subset of $S(\O)$ for every $h$;
\item 
$\{S_h(\O)\}$ is increasing and $\cup_h S_h(\O)=S(\O)$;
\item 
$|D^J\O|(I_L\sm S_h(\O))\to 0$ as $h\to\i$.
\end{itemize}

For every $n,\,h\in \Nat^+$, we estimate
$$ ||DG_n|(\O)-|D\hat G_h|(\O)|\leq |D^J\O|(I_L\sm S_n(\O))+|D^J\O|(I_L\sm S_h(\O)) $$
and hence for every $\e>0$ we can find $n_\e$ and $h_\e$ such that
$||DG_n|(\O)-|D\hat G_h|(\O)| < \e$ if $n>n_\e$ and $h>h_\e$.
Denote as before $\hat\gt_h(s)=\hat G_h(s){\bf e}_1$, and recall that $\gc_n$ is given by \eqref{cn} and $\hat\gc_h$ is defined in a similar way, with $\gt_n$ replaced by $\hat\gt_h$.
We thus can find $\bar n_\e$ and $\bar h_\e$ such that $d(\gc_n,\hat\gc_h)<\e$ if $n>\bar n_\e$ and $h>\bar h_\e$, and hence the sequence $\{\hat\gc_h\}_h$ converges uniformly to the curve $\gc$ previously obtained.

\smallskip 

We now see that a regularization approach leads to the same solution curve $\gc$. 
More precisely, for $\e>0$ small, we can define $\t_\e,\f_\e:\ol I\to\gR$ by the Lipschitz-continuous average functions as in \eqref{approx} and solve the corresponding Cauchy problem \eqref{CP-example}, obtaining the smooth curve $\gc_\e$. 
\bp\label{P-ce} We have 
$\ds\lim_{\e\to 0^+}d(\gc_\e,\gc)=0$.
\ep
\bpf 
We can write the $Sk(3)$-valued $BV$ function $\O$ as a sum $\O=\bar\O+\O^J$, where $\bar\O$ is continuous and $\O^J$ is a purely jump function, in such a way that $\bar D\O=D\bar\O$ and $D^J\O=D\O^J$. For every $n\in\Nat^+$ we let $\O^J_n$ be the purely jump function such that $D\O^J_n=(D^J\O)\pri S_n(\O)$, and we denote by $\t_n$ and $\f_n$ the entries of the matrix valued function $\O_n:=\bar\O+\O^J_n$. 

For $\e>0$ small, we define $\t_{n,\e},\f_{n,\e}:\ol I\to\gR$ by the Lipschitz-continuous average functions as in \eqref{approx}, but with respect to the data $\t_n$ and $\f_n$, and we can solve the corresponding Cauchy problem \eqref{CP-example}, obtaining the smooth curve $\gc_{n,\e}$. 
As in Remark \ref{Rapprox}, since $\O_n$ is continuous outside a finite set, it turns out that $$ \lim_{\e\to 0^+}d(\gc_{n,\e},\gc_n)=0\quad\fa\,n\in\Nat^+\,. $$
On the other hand, for every $s\in I_L$ we can estimate
$$ |\t_\e(s)-\t_{n,\e}(s)|\leq \frac 1{2\e}\int_{s-\e}^{s+\e}|\t(\s)-\t_n(\s)|\,\dd\s\,,
\quad |\f_\e(s)-\f_{n,\e}(s)|\leq \frac 1{2\e}\int_{s-\e}^{s+\e}|\f(\s)-\f_n(\s)|\,\dd\s\,.
$$
Therefore, by the $L^1$ convergences $\t_n\to\t$ and $\f_n\to\f$ we obtain that
$$ \lim_{n\to\i}d(\gc_{n,\e},\gc_\e)=0\quad\textrm{uniformly on }\e>0\,. $$
Writing
$$ d(\gc_\e,\gc)\leq d(\gc_\e,\gc_{n,\e})+d(\gc_{n,\e},\gc_n)+d(\gc_n,\gc)\,, $$
the limit readily follows. 
\epf
\par On account of Proposition \ref{P-ce}, we may thus conclude that the curve $\gc$ given by Proposition \ref{P-cn} corresponds in a weak sense to the (essentially) unique curve whose curvature and torsion are related to the solution of the system
$$
\left\{ \begin{array}{l}
DG= - (G^\top)^{-1} \left(\bar D\O+D^{at}\O\right) \\
G(0)=\bI_3\,.
\end{array} \right.
$$
%
%
\begin{appendices}
\numberwithin{equation}{section}
\section{Proof of the Cayley's theorem}
\bpf[Proof of Proposition \ref{P-Cayley}:] Using that $(M^\top)^{-1}=(M^{-1})^\top$ for invertible matrices $M\in\gR^{N\tim N}$, and letting $\bI=\bI_N$, we have
\[
B^\top = (A-\bI)^\top((A+\bI)^{-1})^{\top}=(A^\top-\bI)(A^\top + \bI)^{-1}
\]
whence
\[
B+B^\top=(A+\bI)^{-1}(A-\bI)+(A^\top-\bI)(A^\top + \bI)^{-1} .
\]
We thus can write
\[
(B+B^\top)=(A+\bI)^{-1}\,\left[(A-\bI)+(A+\bI)(A^\top-\bI)(A^\top + \bI)^{-1} \right]
\]
where
\[
(A+\bI)(A^\top-\bI)=A^\top-A =(\bI-A)(A^T+\bI)
\]
so that we obtain
\[
\ba{rl}
(B+B^\top)= & (A+\bI)^{-1}\,\left[(A-\bI)+(\bI-A)(A^\top+\bI)(A^\top + \bI)^{-1} \right] \\
= &
(A+\bI)^{-1}\,\left[(A-\bI)+(\bI-A)\right]=0\,,
\ea \]
as required. \epf
%
%
\section{Smooth solutions to the Frenet-Serret system}\label{App:B}
In the smooth case when $\O \in C^1(I_L,Sk(3))$, in general one cannot find the explicit formula of smooth solutions $G\in C^1(I_L,SO(3))$ to equation
\beq\label{eq-smooth}
G'(s)=-G(s)\O'(s)\,.
\eeq

In fact, any map $G\in C^1(I_L,SO(3))$ is given by $G(s)=\exp( X(s))$ for some map $X\in C^1(I_L,Sk(3))$, and
the general derivative formula gives
\beq\label{deriv}
\frac d{ds}\exp(X(s))=\exp(X(s))\,B(s)\,,\quad B(s):=\int_0^1 \exp(-\s X(s))X'(s)\,\exp(\s X(s))\,\dd \s\,.
\eeq
Therefore, the required solution $G(s)$ corresponds to a function $B(s)$ satisfying equation
\beq\label{B-Om}
B(s)=-\O'(s)\quad\fa\,s\in I_L\,.
\eeq

If the product $X(s_1)X(s_2)$ is commutative for every $s_1,s_2 \in I_L$, we readily obtain that $B(s)=X'(s)$, whence any solution to eq. \eqref{eq-smooth} is given by $G(s)=\exp(-\O(s)+c)$ for some constant $c\in\gR$.
In general, we have $B(s)\neq X'(s)$, and one cannot find the explicit formula of solutions to eq. \eqref{eq-smooth}. However, something more can be said.

In fact, we can write $X(s)$ in terms of a triplet of smooth functions $f,g,h:I_L\to\gR$ through formula $X(s)= \r(s)\bJ(s)$, where $\r(s)=\left\{ f(s)^2+g(s)^2+h(s)^2 \right\}^{1/2}$ and
$$ \bJ(s) := \left(\begin{array}{ccc}
0 & \a(s) & \g(s) \\ -\a(s) & 0 & \be(s) \\ -\g(s)  & -\be(s) & 0  \\
 \end{array} \right) \in Sk(3)\,, $$
with $\a(s)=f(s)/\r(s)$, $\be(s)=g(s)/\r(s)$, $\g(s)=h(s)/\r(s)$.
On account of eq. \eqref{matrix-J}, we have
$$ \bJ(s)=-\be(s)\,\bJ_1+\g(s)\,\bJ_2-\a(s)\,\bJ_3\,,\quad s \in I_L\,.
$$
Therefore, Rodrigues formula \eqref{Rodrigues} gives for every $\l\in\gR$
$$ \exp(\lambda X)=\left(\begin{array}{ccc}
c+(1-c)\be^2 & s\,\a-(1-c)\be\g  & s\,\g+(1-c)\a\be \\ -s\,\a-(1-c)\be\g & c+(1-c)\g^2 & s\,\be-(1-c)\a\g \\ -s\,\g+(1-c)\a\be  & -s\,\be-(1-c)\a\g & c+(1-c)\a^2  \\
 \end{array} \right)
$$
where $s=\sin\l$, $c=\cos\l$. As a consequence, the matrix-valued function $B(s)\in C^1(\ol I_L,Sk(3))$ in the derivative formula \eqref{deriv} is identified by the entries $B_{12}(s)$, $B_{13}(s)$, $B_{23}(s)$ that are given by
\beq\label{B}
\left\{ \ba{ll}
B_{12}(s) = & \ds \sin\r(s)\frac \dd{\dd s}\frac{f(s)}{\r(s)}+\frac{f(s)}{\r(s)}\,\r'(s)+\frac{1-\cos\r(s)}{\r(s)^2}\,\left( g'(s)h(s)-g(s)h'(s)\right)
\\
B_{13}(s) = & \ds \sin\r(s)\frac \dd{\dd s}\frac{h(s)}{\r(s)}+\frac{h(s)}{\r(s)}\,\r'(s)+\frac{1-\cos\r(s)}{\r(s)^2}\,
\left( g'(s)f(s)-g(s)f'(s)\right)
\\
B_{23}(s) = & \ds \sin\r(s)\frac \dd{\dd s}\frac{g(s)}{\r(s)}+\frac{g(s)}{\r(s)}\,\r'(s) +\frac{1-\cos\r(s)}{\r(s)^2}\,
\left( f'(s)h(s)-f(s)h'(s)\right)\,.
%
\ea \right.
\eeq

Notice that if $f(s)=ks$, $g(s)=\tau s$ and $h(s)=0$, with $k>0$ and $\tau\in\gR$, we obtain $B_{12}(s)=k$, $B_{13}(s)=0$, $B_{23}(s)=\tau$ for every $s$. However, given $\t'(s)=k(s)>0$ and $\f'(s)=\tau(s)$ smooth, in general one cannot find the explicit formulas of the function $f$, $g$, and $h$ in \eqref{B} in such a way that eq. \eqref{B-Om} holds, i.e.
$$
\left\{ \ba{ll}
B_{12}(s) = & k(s)
 \\
B_{13}(s) = & 0
 \\
B_{23}(s) = & \tau(s)\,,%
\ea \right. \quad s\in I_L\,.
$$
\end{appendices}

\end{document}